\documentclass[12pt,reqno]{amsart}

\title[]{Long Time Behavior of Solutions of an Electroconvection Model in $\R^2$}

\author{Elie Abdo}
\address{Department of Mathematics, Temple University, Philadelphia, PA 19122}
\email{abdo@temple.edu}
\author{Mihaela Ignatova}
\address{Department of Mathematics, Temple University, Philadelphia, PA 19122}
\email{ignatova@temple.edu}

\usepackage[margin=1in]{geometry}
\usepackage{amsmath, amsthm, amssymb}
\usepackage{times}
\usepackage{color}
\usepackage{hyperref}
\newcommand{\pa}{\partial}
\newcommand{\la}{\label}
\newcommand{\fr}{\frac}
\newcommand{\na}{\nabla}
\newcommand{\be}{\begin{equation}}
\newcommand{\ee}{\end{equation}}
\newcommand{\ba}{\begin{array}{l}}
\newcommand{\ea}{\end{array}}
\newcommand{\Rr}{{\mathbb R}}

\newcommand{\beg}{\begin}

\renewcommand{\l}{\Lambda}

\usepackage{MnSymbol,wasysym}

\newcommand{\R}{\mathbb R}

\def\RR{{\mathbb R}}

\def\PP{\mathbb P}

\date{today}
\begin{document}
\begin{abstract} 
We consider a two dimensional electroconvection model which consists of a nonlinear and nonlocal system coupling the evolutions of a charge distribution and a fluid. We show that the solutions decay in time in $L^2(\Rr^2)$ at the same sharp rate as the linear uncoupled system. This is achieved by proving that the difference between the nonlinear and linear evolution decays at a faster rate than the linear evolution. In order to prove the sharp $L^2$ decay we establish bounds for decay in $H^2(\Rr^2)$ and a logarithmic growth in time of a quadratic moment of the charge density.
 \end{abstract} 

\maketitle
\section{Introduction}\la{intro}
We consider the electroconvection model 
\be \la{intro1}
\pa_t q + u \cdot \na q + \l q = 0,
\ee
\be  \la{intro2}
\pa_t u + u \cdot \na u + \na p - \Delta u = -qRq,
\ee 
\be \la{intro3}
\na \cdot u = 0
\ee
in $\R^2$ describing the evolution of a surface charge density $q$ in a two-dimensional incompressible fluid flowing with a velocity $u$ and a pressure $p$. Here $\l = (-\Delta)^{\fr{1}{2}}$ is the square root of the two-dimensional Laplacian, and $R = \na \l^{-1}$ is the two-dimensional Riesz transform. The initial data 
\be 
q(\cdot, 0) = q_0
\la{qinit}
\ee
and
\be
u(\cdot, 0) = u_0
\la {uinit}
\ee
are assumed to be regular enough and have good decay properties.
The model is motivated by physical and numerical studies of electroconvection
\cite{DDMB,DMT,DDMT}. The nonlocal aspect of the evolution of the charge density
and the nonlocal forcing on the Navier-Stokes equations in the model are due to the fact that the fluid and charges are confined to a thin two dimensional film. The global well-posedness of the system in bounded domains was obtained in \cite{ceiv} using commutator estimates and nonlocal nonlinear analysis.
In \cite{AI}, we investigated the long time dynamics of the model in two dimensions, with periodic boundary conditions and with applied voltage. When the fluid is forced by time-independent smooth mean zero body forces, we proved that the model \eqref{intro1}--\eqref{qinit} has a finite dimensional global attractor. In the absence of body forces, the charge density $q$ converges exponentially in time to a unique limit due to the applied voltage, and the velocity $u$ converges exponentially in time to zero. The rate of exponential decay depends  on the periodic boundary conditions.

In this paper, we consider the time asymptotic behavior of solutions of \eqref{intro1}--\eqref{qinit} in $\RR^2$, and adapt the Fourier splitting method \cite{S,S1} of Schonbek to the present system. The method was initially used in \cite{S} to prove decay of Leray weak solutions \cite{L} of Navier-Stokes equations and to further decay studies for Navier-Stokes equations \cite{BN,KM,S1,SW,W1} and 
many other partial differential equations (see for instance \cite{BS,CW,DL,NS,Z,ZL}). Different approaches were employed as well to investigate the time decay \cite{OT} and space-time decay  \cite{AGSS, KU1, KU,T} of higher-order derivatives of solutions to Navier-Stokes equations.

The electroconvection model \eqref{intro1}--\eqref{qinit} couples Navier-Stokes equations to a scalar equation for a surface charge density $q$, evolving via advection by $u$ and diffusion by $\Lambda $.  We obtain in Theorem~\ref{decay} of section \ref{s2} the long time $L^2$ decay of the type
\[
\|q\|_{L^2} = O(t^{-1})
\]
and
\[
\|u\|_{L^2}  = O(t^{-\fr{1}{2}}).
\]
This rate of decay is sharp for the linear uncoupled system if the initial data have non vanishing finite $L^1$ norms, because functions of the form  $Q(t) = e^{-t\l^\alpha}q_0$  obey
\[
\lim_{t\to\infty} t^{\frac{n}{\alpha}}\|Q(t)\|_{L^2(\Rr^n)}^2 = C_{n,\alpha} \left(\int_{\Rr^n} q_0 dx\right)^2
\]
for any $\alpha>0$ and $n\ge 1$.
The fact that such a decay is sharp for the nonlinear evolution as well is a consequence of Theorem~\ref{nonlin} of section \ref{s3} where we prove that $u- U$  with $U(t) = e^{t\Delta}u_0$ and $q- Q$ with $Q(t) = e^{-t\l}q_0$ decay faster in $L^2$ than $u$ and $q$, respectively. Similar results were proved for critical SQG in \cite{CW}. The critical SQG velocity $u= R^{\perp}q$ decays in $L^2$ like $q$, that is at the rate $t^{-1}$, which helps lower the size of the nonlinear term $u\cdot\na q$ in that equation. In our case, the velocity has slower decay in $L^2$ due to the Navier-Stokes equation, namely of the order $t^{-\fr{1}{2}}$, and the nonlinear term is larger. The influence of the charge density $q$ is felt by the Navier-Stokes velocity via the electric force $-qRq$. In order to obtain a key fast enough decay at low wave numbers for the difference 
$v= u- U$, we need to control a moment of $q$,  $\int_{\Rr^2}|x|^2|q(x,t)|^2dx = M^2(t)$, in view of the inequality
\[
|\widehat{qRq}(\xi)| \le C |\xi| \|q\|_{L^2}M(t)
\]
(see Lemma \ref{lem3}).  We prove that 
\[
M(t) = O(\sqrt{\log t})
\]
for long time, by analyzing the evolution of the quantity $a(x)q(x,t)$ with $a(x)=\sqrt{|x|^2 +1}$. This analysis uses the boundedness of the commutator between $\Lambda$ and multiplication by $a(x)$, which we establish in Lemma \ref{lem1}. In addition, in order to achieve the necessary sharp $L^2$ bounds we have to obtain bounds for the decay of higher norms of both $u$ and $q$.  For instance,  $H^1$
 norms of $q$ are of the order 
\[
\|\na q\|_{L^2} = O(t^{-1}).
\]
These bounds are obtained by somewhat involved nonlinear and nonlocal analysis, and they are no longer sharp compared to the generic $t^{-2}$ linear behavior.

The paper is organized as follows. In section \ref{s2}, we study the asymptotic behavior of solutions to the electroconvection model \eqref{intro1}--\eqref{qinit}: we prove that the $L^2$ norm of the surface charge density $q$ decays in time to zero with a rate of order $t^{-1}$ whereas the velocity $u$ decays in time to zero with a rate of order  $t^{-\fr{1}{2}}$. We also investigate the rate of decay of their higher-order derivatives, and we obtain decaying-in-time bounds in H\"older spaces $C^{0, \fr{1}{2}}$. In section \ref{s3}, we prove that the differences $q-Q$ and $u - U$ decay to zero in $L^2$ faster than $q$ and $u$, with rates of order $t^{-1- \fr{3}{4}}$ and $t^{-\fr{3}{4}}$, respectively. In the Appendix, we present results on the existence and uniqueness of solutions to \eqref{intro1}--\eqref{qinit}, based on the Banach fixed point theorem, the Aubin-Lions lemma and commutator estimates.

\section{Long Time Behavior of Solutions} \la{s2}

In this section, we consider the long-time behavior of solutions of the electroconvection model described by \eqref{intro1}--\eqref{qinit}. We show that the charge density $q$ and the velocity $u$ converge to $0$ in the $H^2$ norm, and we investigate the rate of convergence.  

For a function $f \in L^1(\RR^2)$, we denote its Fourier transform by 
\be 
\widehat{f} (\xi) = \int_{\RR^2} f(x) e^{- i \xi \cdot x } dx.
\ee

\beg{thm} \la{decay} Let $u_0 \in H^1 \cap L^1$ be divergence-free and $q_0 \in L^4 \cap L^1$. There exist positive constants $\Gamma_0$ and $\Gamma_0'$ depending only on the initial data and some universal constants such that the unique global-in-time solution $(q,u)$ of \eqref{intro1}--\eqref{qinit} obeys
\be \la{ltime1}
\|q(t)\|_{L^2}^2 \le \fr{\Gamma_0}{(t + 1)^{2}}
\ee and 
\be \la{ltime2}
\|u(t)\|_{L^2}^2 \le \fr{\Gamma_0'}{t+1}
\ee
for all $t \ge 0$. 
\end{thm}

\textbf{Proof:} The proof is divided into several steps. 

\textit{\textbf{Step 1} (Basic energy estimates}).
We take the $L^2$ inner product of equation \eqref{intro1} with $\l^{-1} q$ and the $L^2$ inner product of equation \eqref{intro2} with $u$. Then we add the resulting energy equalities. Integrating by parts, we have the cancellations 
\be 
(u \cdot \na u, u)_{L^2} = (\na p, u)_{L^2} = 0
\ee and 
\beg{align}
(u \cdot \na q, \l^{-1} q)_{L^2} + (qRq, u)_{L^2} 
&= - (u \cdot \na \l^{-1} q, q)_{L^2} + (qRq, u)_{L^2} \nonumber
\\&= - (u \cdot Rq, q)_{L^2} + (qRq, u)_{L^2}
= 0
\end{align} due to the divergence-free condition \eqref{intro3}.
Thus, we obtain 
\be \la{weaknormb}
\fr{1}{2} \fr{d}{dt} \left(\|\l^{-\fr{1}{2}}q\|_{L^2}^2 + \|u\|_{L^2}^2  \right)
+ \|q\|_{L^2}^2 + \|\na u\|_{L^2}^2
= 0.
\ee We integrate in time from $0$ to $t$ and we take the supremum over all positive times $t \ge 0$. We get 
\be  \la{step11}
\sup\limits_{t \ge 0} \left\{\|\l^{-\fr{1}{2}} q(t)\|_{L^2}^2 + \|u(t)\|_{L^2}^2
+ \int_{0}^{t} 2\left(\|q(s)\|_{L^2}^2 + \|\na u(s)\|_{L^2}^2 \right) ds  \right\}
= \|\l^{-\fr{1}{2}} q_0\|_{L^2}^2 + \|u_0\|_{L^2}^2
\ee ending the proof of Step 1.

\textit{\textbf{Step 2} (Pointwise bounds for the Fourier transform of the charge density $q$).}
The Fourier transform of $q$ evolves according to 
\be 
\pa_t \widehat{q} (\xi, t) + (\widehat{u \cdot \na q})(\xi, t) + \widehat{\l q} (\xi, t) = 0.
\ee The fractional Laplacian $\l$ is a Fourier multiplier with symbol $|\xi|$, hence 
\be 
\pa_t \widehat{q} + |\xi| \widehat{q} = - \widehat{u \cdot \na q}.
\ee We estimate the Fourier transform of the nonlinear term 
\be 
|\widehat{u \cdot \na q}| = |\widehat{\na \cdot (uq)}| \le C|\xi| \|u\|_{L^2} \|q\|_{L^2}
\ee using the divergence-free condition \eqref{intro3}, the boundedness of the Fourier transform of a function by its $L^1$ norm, and the Cauchy-Schwarz inequality. 
This yields the differential inequality 
\be 
\pa_t \widehat{q} + |\xi| \widehat{q}  \le C|\xi| \|u\|_{L^2} \|q\|_{L^2}.
\ee We multiply both sides by the integrating factor $e^{|\xi|t}$ and integrate in time from $0$ to $t$. We obtain the bound
\be 
|\widehat{q} (\xi, t)| \le |\widehat{q}_0 (\xi)| + C|\xi| \int_{0}^{t} \|u(s)\|_{L^2}\|q(s)\|_{L^2} ds. 
\ee As a consequence of Step 1 and the Cauchy-Schwarz inequality, we get the pointwise bound
\be \la{EQ1}
|\widehat{q} (\xi, t)| \le \|q_0\|_{L^1} + C_0|\xi|\sqrt{t} 
\ee where $C_0$ is a time-independent constant depending only on $\|u_0\|_{L^2}$ and $\|\l^{-\fr{1}{2}}q_0\|_{L^2}$. This finishes the proof of Step 2.

\textit{\textbf{Step 3} (Decaying bound for the $L^2$ norm of the charge density).}
The $L^2$ norm of $q$ evolves according to 
\be \la{step22}
\fr{1}{2} \fr{d}{dt} \|q\|_{L^2}^2 + \|\l^{\fr{1}{2}} q\|_{L^2}^2 = 0.
\ee In view of Parseval's identity and the fact that $\l^{\fr{1}{2}}$ is a Fourier multiplier with symbol $|\xi|^{\fr{1}{2}}$, we have 
\be 
\|\l^{\fr{1}{2}} q\|_{L^2}^2 = \|\widehat{\l^{\fr{1}{2}}q}\|_{L^2}^2 = \int_{\R^2} |\xi| |\widehat{q}(\xi,t)|^2 d\xi. 
\ee We bound the dissipation from below 
\be 
\int_{\R^2} |\xi| |\widehat{q}(\xi,t)|^2 d\xi 
\ge \int_{|\xi| > \rho(t)} |\xi| |\widehat{q}(\xi,t)|^2 d\xi 
\ee where $\rho(t)$ is the function defined on $[0, \infty)$ by 
\be 
\rho(t) = \frac{r}{2(t+1)}
\ee for some positive constant $r$ to be determined later. We note that 
\beg{align} 
\int_{|\xi| > \rho (t)} |\xi| |\widehat{q}(\xi,t)|^2 d\xi 
&\ge \rho(t) \int_{|\xi| > \rho(t)} |\widehat{q}(\xi, t)|^2 d\xi \nonumber
\\&= \rho(t) \int_{\R^2} |\widehat{q}(\xi, t)|^2 d\xi - \rho(t) \int_{|\xi| \le \rho(t)} |\widehat{q}(\xi, t)|^2 d\xi \nonumber
\\&= \rho(t) \|q\|_{L^2}^2 - \rho(t) \int_{|\xi| \le \rho(t)} |\widehat{q}(\xi, t)|^2 d\xi
\end{align} where we used Parseval's identity. 
Consequently, we obtain the energy inequality
\be 
\fr{d}{dt} \|q\|_{L^2}^2 + 2\rho(t) \|q\|_{L^2}^2 \le 2\rho(t) \int_{|\xi| \le \rho(t)} |\widehat{q}(\xi, t)|^2 d\xi.
\ee
By the pointwise bound \eqref{EQ1} and Fubini's theorem for spherical coordinates, we estimate
\beg{align}
\int_{|\xi| \le \rho(t)} |\widehat{q}(\xi, t)|^2 d\xi
&\le \int_{|\xi| \le \rho(t)} \left(\|q_0\|_{L^1} + C_0|\xi|\sqrt{t} \right)^2 d\xi
= C\int_{0}^{\rho(t)} r\left(\|q_0\|_{L^1} + C_0 r\sqrt{t} \right)^2 dr \nonumber
\\&\le C\int_{0}^{\rho(t)} r\left(\|q_0\|_{L^1}^2 + C_0^2 r^2 t  \right) dr
\le \Gamma_1 \left(\rho(t)^2 + t\rho(t)^4 \right) 
\end{align} where $\Gamma_1$ depends only on the initial data. 
We obtain 
\be 
\fr{d}{dt} \|q\|_{L^2}^2 + 2\rho(t) \|q\|_{L^2}^2 \le 2\Gamma_1 (\rho(t)^3 + t\rho(t)^5) 
\ee for all $t \ge 0$. We multiply both sides of the inequality by the integrating factor 
\be 
e^{2\int_{0}^{t} \rho(s) ds} = e^{r \int_{0}^{t} \fr{1}{s+1} ds} = e^{r \ln(t+1)} = (t+1)^{r}
\ee and then we integrate in time from $0$ to $t$. We get
\beg{align} 
\|q(t)\|_{L^2}^2 
&\le \fr{\|q_0\|_{L^2}^2}{(t+1)^{r}} + \fr{\Gamma_2}{(t+1)^{r}} \int_{0}^{t} \left(\fr{1}{(s+1)^3} + \fr{1}{(s+1)^4}\right) (s+1)^{r} ds \nonumber
\\&\le \fr{\|q_0\|_{L^2}^2}{(t+1)^r} + \fr{\Gamma_2}{(t+1)^{r}} \left(\fr{(t+1)^{r-2}}{r-2} - \fr{1}{r-2} + \fr{(t+1)^{r-3}}{r-3} - \fr{1}{r-3} \right) \nonumber
\\&\le \fr{\|q_0\|_{L^2}^2}{(t+1)^{r}} + \fr{\Gamma_2}{(r-2)(t+1)^{2}} + \fr{\Gamma_2}{(r-3)(t+1)^{3}} \la{step33}
\end{align} for any $r > 3$. Here $\Gamma_2$ depends on $r$ and the initial data. We choose $r=4$ and we obtain the bound
\be \la{EQ3}
\|q\|_{L^2}^2 \le \fr{\Gamma_3}{(t+1)^2}
\ee where $\Gamma_3$ is a positive constant depending only on the initial data. This completes the proof of \eqref{ltime1} and Step~3. 

\textit{\textbf{Step 4} (Pointwise bounds for the Fourier transform of the velocity $u$).}
Applying the Leray Projector $\PP$ to equation \eqref{intro2}, we have 
\be 
\pa_t u + \PP(u \cdot \na u) - \Delta u = -\PP(qRq),
\ee where we used the incompressibility condition \eqref{intro3} and the fact that $\PP$ and $-\Delta$ are Fourier multipliers so they commute. 
Hence the Fourier transform of $u$ obeys 
\be 
\pa_t \widehat{u} + \widehat{ \PP(u \cdot \na u)} - \widehat{\Delta u} = -\widehat{\PP(qRq)}.
\ee 
We estimate 
\be 
|\widehat{ \PP(u \cdot \na u)}(\xi, t)| \le C|\xi||\widehat{u}(\xi,t)|^2 \le C|\xi|\|u(t)\|_{L^2}^2 
\ee and 
\be 
|\widehat{\PP(qRq)}(\xi,t)| \le C\|(qRq)(t)\|_{L^1} \le C\|q(t)\|_{L^2}^2 
\ee in view of the boundedness of the Riesz transforms on $L^2(\RR^2)$.
We obtain 
\be 
\pa_t \widehat{u} + |\xi|^2 \widehat{u} \le C|\xi|\|u\|_{L^2}^2 + C\|q\|_{L^2}^2
\ee and hence
\be 
|\widehat{u}(\xi,t)| \le \|u_0\|_{L^1} + C|\xi| \int_{0}^{t} \|u(s)\|_{L^2}^2 ds + C\int_{0}^{t}\|q(s)\|_{L^2}^2 ds
\ee for all $\xi \in \RR^2$ and $t \ge 0$. 
In view of the bound \eqref{step11}, we get 
\be \la{EQ2}
|\widehat{u}(\xi,t)| \le \Gamma_4  + C|\xi| \int_{0}^{t} \|u(s)\|_{L^2}^2 ds 
\ee where $\Gamma_4$ is a positive constant depending only on the initial data. This completes the proof of Step~4.

\textit{\textbf{Step 5} (Decaying bounds for the $L^4$ norm of $q$).} The evolution of the $L^4$ norm of $q$ is described by the energy equality
\be 
\fr{1}{4} \fr{d}{dt} \|q\|_{L^4}^4 + \int_{\R^2} q^3 \l q dx = 0.
\ee In view of the C\'ordoba-C\'ordoba inequality \cite{CC}, the dissipation is bounded from below 
\be 
\int_{\R^2} q^3 \l q dx 
\ge \fr{1}{2} \|\l^{\fr{1}{2}}(q^2)\|_{L^2}^2 
\ee and thus
\be 
\int_{\R^2} q^3 \l q dx  \ge c\|q\|_{L^8}^4 
\ee due to Gagliardo-Nirenberg inequalities. Using interpolation inequalities in $L^p$ spaces and the uniform boundedness of the $L^2$ norm of the charge density $q$ by $\|q_0\|_{L^2}$, we have the bound
\be 
\|q\|_{L^4} \le \|q\|_{L^2}^{\fr{1}{3}}  \|q\|_{L^8}^{\fr{2}{3}}\le \|q_0\|_{L^2}^{\fr{1}{3}} \|q\|_{L^8}^{\fr{2}{3}}
\ee from which we conclude that
\be 
\int_{\R^2} q^3 \l q dx \ge C \|q_0\|_{L^2}^{-2} \|q\|_{L^4}^6
\ee and hence
\be 
\fr{d}{dt} \|q\|_{L^4}^4 + \fr{C}{\|q_0\|_{L^2}^2} \|q\|_{L^4}^6 \le 0.
\ee Letting $y = \|q\|_{L^4}^4$, we obtain the Bernouilli ordinary differential inequality 
\be 
\fr{dy}{dt}  + \fr{C}{\|q_0\|_{L^2}^2} y^{\fr{3}{2}} \le 0.
\ee We apply a change of variable given by $u = y^{-\fr{1}{2}}$ and we get
\be 
\fr{-2}{u^3} \fr{du}{dt}  + \fr{C}{\|q_0\|_{L^2}^2} \fr{1}{u^3} \le 0
\ee so
\be 
\fr{du}{dt} \ge \fr{C}{2\|q_0\|_{L^2}^2}.
\ee Integrating in time from $0$ to $t$, we arrive at the bound
\be 
\|q\|_{L^4}^{-2} \ge \|q_0\|_{L^4}^{-2} + \fr{C}{\|q_0\|_{L^2}^2}t \ge \Gamma_5(1+t) 
\ee where $\Gamma_5$ is a constant depending on the initial data. Consequently, we obtain 
\be \la{Bernouilli}
\|q\|_{L^4} \le \fr{1}{\sqrt{\Gamma_5}} \fr{1}{(1+t)^{\fr{1}{2}}}
\ee for all $t \ge 0$.

\textit{\textbf{Step 6} (Decaying bound for the $L^2$ norm of the velocity).}
The $L^2$ norm of the velocity evolves according to 
\be 
\fr{1}{2} \fr{d}{dt} \|u\|_{L^2}^2 + \|\na u\|_{L^2}^2 = -\int_{\RR^2} qRq \cdot u dx.
\ee
In view of H\"older's inequality, the boundedness of the Riesz transforms on $L^4(\RR^2)$, and Ladyzhenskaya's interpolation inequality, we bound
\be 
\left|\int_{\RR^2} qRq \cdot u \right| \le C\|q\|_{L^2}\|q\|_{L^4}\|u\|_{L^2}^{\fr{1}{2}}\|\na u\|_{L^2}^{\fr{1}{2}}
\ee yielding 
\be 
 \fr{d}{dt} \|u\|_{L^2}^2 + \|\na u\|_{L^2}^2 \le C\|q\|_{L^2}^{\fr{4}{3}} \|q\|_{L^4}^{\fr{4}{3}} \|u\|_{L^2}^{\fr{2}{3}}.
\ee By Parseval's identity, we have
\be 
\fr{d}{dt} \|u\|_{L^2}^2 + \int_{\R^2} |\xi|^2 |\widehat{u}(\xi,t)|^2 d\xi \le C\|q\|_{L^2}^{\fr{4}{3}} \|q\|_{L^4}^{\fr{4}{3}} \|u\|_{L^2}^{\fr{2}{3}}.
\ee For a positive function $\rho_1(t)$ continuous on $[0, \infty)$, we have  
\beg{align}
\int_{\R^2} |\xi|^2 |\widehat{u}(\xi,t)|^2 d\xi
&\ge \int_{|\xi| > \rho_1(t)} |\xi|^2 |\widehat{u}(\xi,t)|^2 d\xi
\ge \rho_1(t)^2 \int_{|\xi| > \rho_1(t)} |\widehat{u}(\xi,t)|^2 d\xi \nonumber
\\&\ge \rho_1(t)^2 \left(\int_{\R^2} |\widehat{u}(\xi,t)|^2 d\xi - \int_{|\xi| \le \rho_1(t)}  |\widehat{u}(\xi,t)|^2 d\xi \right) \nonumber
\\&= \rho_1(t)^2 \|u\|_{L^2}^2 - \rho_1(t)^2 \int_{|\xi| \le \rho_1(t)}  |\widehat{u}(\xi,t)|^2 d\xi.
\end{align} 
Consequently, we obtain the energy inequality 
\be 
\fr{d}{dt} \|u\|_{L^2}^2 + \rho_1(t)^2 \|u\|_{L^2}^2 \le C\|q\|_{L^2}^{\fr{4}{3}} \|q\|_{L^4}^{\fr{4}{3}} \|u\|_{L^2}^{\fr{2}{3}} + \rho_1(t)^2\int_{|\xi| \le \rho_1(t)}  |\widehat{u}(\xi,t)|^2 d\xi.
\ee
Using \eqref{EQ2}, we have 
\beg{align} 
&\int_{|\xi| \le \rho_1(t)}  |\widehat{u}(\xi,t)|^2 d\xi
\le C\int_{0}^{\rho_1(t)} r \left(\Gamma_4^2 + C r^2 \left\{\int_{0}^{t} \|u(s)\|_{L^2}^2 ds \right\}^2 \right) dr \nonumber
\\&\le \Gamma_6 \rho_1(t)^2 +  C\rho_1(t)^4 \left(\int_{0}^{t} \|u(s)\|_{L^2}^2 ds  \right)^2 
\end{align} and thus
\beg{align} 
\fr{d}{dt} \|u\|_{L^2}^2 + \rho_1(t)^2 \|u\|_{L^2}^2 \le  \Gamma_6 \rho_1(t)^4 + C\rho_1(t)^6  \left(\int_{0}^{t} \|u(s)\|_{L^2}^2 ds \right)^2 + C\|q\|_{L^2}^{\fr{4}{3}} \|q\|_{L^4}^{\fr{4}{3}} \|u\|_{L^2}^{\fr{2}{3}}.
\end{align}
Multiplying by the integrating factor $e^{\int_{0}^{t} \rho_1(s)^2 ds}$, and integrating in time from $0$ to $t$, we obtain
\beg{align} 
\|u(t)\|_{L^2}^2 
&\le \fr{\|u_0\|_{L^2}^2}{e^{\int_{0}^{t} \rho_1(s)^2 ds}}
+ \fr{\Gamma_{6}}{e^{\int_{0}^{t} \rho_1(s)^2 ds}} \int_{0}^{t} e^{\int_{0}^{s} \rho_1(\tau)^2 d\tau} \rho_1(s)^4 ds  \nonumber
\\&+ \fr{C}{e^{\int_{0}^{t} \rho_1(s)^2 ds}} \int_{0}^{t} \left(e^{\int_{0}^{s} \rho_1(\tau)^2 d\tau}\rho_1(s)^6\right) \left(\int_{0}^{s} \|u(\tau)\|_{L^2} ^2 d\tau \right)^2 ds  \nonumber
\\&+ \fr{C}{e^{\int_{0}^{t} \rho_1(s)^2 ds}} \int_{0}^{t} \|q\|_{L^2}^{\fr{4}{3}} \|q\|_{L^4}^{\fr{4}{3}} \|u\|_{L^2}^{\fr{2}{3}} e^{\int_{0}^{s} \rho_1(\tau)^2 d\tau} ds.  
\end{align} 
In view of 
\eqref{step11}, \eqref{EQ3} and \eqref{Bernouilli}, we estimate 
\be 
 \int_{0}^{t} \|q\|_{L^2}^{\fr{4}{3}} \|q\|_{L^4}^{\fr{4}{3}} \|u\|_{L^2}^{\fr{2}{3}} e^{\int_{0}^{s} \rho_1(\tau)^2 d\tau} ds
\le \Gamma_7 \int_{0}^{t} \fr{e^{\int_{0}^{s} \rho_1(\tau)^2 d\tau}}{(s+1)^{2}} ds 
\ee for any $t \ge 0$, and so
\beg{align} \la{abstractbound}
\|u(t)\|_{L^2}^2 
&\le \fr{\|u_0\|_{L^2}^2}{e^{\int_{0}^{t} \rho_1(s)^2 ds}}
+ \fr{\Gamma_{6}}{e^{\int_{0}^{t} \rho_1(s)^2 ds}} \int_{0}^{t} e^{\int_{0}^{s} \rho_1(\tau)^2 d\tau} \rho_1(s)^4 ds  \nonumber
\\&+ \fr{C}{e^{\int_{0}^{t} \rho_1(s)^2 ds}} \int_{0}^{t} \left(e^{\int_{0}^{s} \rho_1(\tau)^2 d\tau}\rho_1(s)^6\right) \left(\int_{0}^{s} \|u(\tau)\|_{L^2} ^2 d\tau \right)^2 ds   \nonumber
\\&+  \fr{\Gamma_7}{e^{\int_{0}^{t} \rho_1(s)^2 ds}} \int_{0}^{t} \fr{e^{\int_{0}^{s} \rho_1(\tau)^2 d\tau}}{(s+1)^{2}} ds 
\end{align} for any $t \ge 0$.

In order to obtain the sharp decaying bound for the velocity $u$, we need the following three sub-steps: 

\textit{\textbf{Step 6.1} (Logarithmic decaying bound for the $L^2$ norm of the velocity).} We take $\rho_1(t) = (e + t)^{-\fr{1}{2}} \left[\ln (e + t)\right]^{-\fr{1}{2}}$. In this case, the integrating factor is given by 
\be 
e^{\int_{0}^{t} \rho_1(s)^2 ds} = e^{\int_{0}^{t} \fr{1}{(e+s)\ln(e+s)} ds} = e^{\ln \left[\ln (e+t) \right]} = \ln (e+t)
\ee and so \eqref{abstractbound} becomes
\beg{align}
&\|u(t)\|_{L^2}^2 
\le \fr{\|u_0\|_{L^2}^2}{\ln (e+t)}
+ \fr{\Gamma_{6}}{\ln (e +t)} \int_{0}^{t} \fr{1}{(e+s)^2 \ln (e+s)} ds  \nonumber
\\&+ \fr{C \|u_0\|_{L^2}^2}{\ln (e +t)} \int_{0}^{t} \fr{s^2}{(e+s)^3 \left[\ln (e+s)\right]^2} ds  
+  \fr{\Gamma_7}{\ln (e+t)} \int_{0}^{t} \fr{\ln (e+s)}{(s+1)^{2}} ds. 
\end{align} in view of the uniform boundedness of $\|u\|_{L^2}$ by $\|u_0\|_{L^2}$. We note that 
\be 
\int_{0}^{t} \fr{s^2}{(e+s)^3 \left[\ln (e+s)\right]^2} ds 
\le \int_{0}^{t} \fr{1}{(e+s) \left[\ln (e+s) \right]^2} ds
= 1 - \fr{1}{\ln (e+t)} \le 1
\ee for any $t \ge 0$. Therefore, 
\be \la{EQ4}
\|u(t)\|_{L^2}^2 \le \fr{\Gamma_8}{\ln (e+t)}
\ee for all $t \ge 0$, where $\Gamma_8$ is a constant depending only on the initial data. 

\textit{\textbf{Step 6.2} (Almost sharp decaying bound for the $L^2$ norm of the velocity).} In order to improve the logarithmic decay \eqref{EQ4}, we take $\rho_1(t) = r^{\fr{1}{2}} (t+1)^{-\fr{1}{2}}$ for some $r$ to be chosen later. In this case, the integrating factor is given by
\be 
e^{\int_{0}^{t} \rho_1(s)^2 ds} = e^{r \int_{0}^{t} \fr{1}{(s+1)} ds} = e^{r \ln (t+1)} = (t+1)^r
\ee and so \eqref{abstractbound} becomes
\be \la{abstractbound1}
\|u(t)\|_{L^2}^2 
\le \fr{\|u_0\|_{L^2}^2}{(t+1)^r}
+ \fr{\Gamma_{9}}{(t+1)^r} \int_{0}^{t} \fr{(s+1)^r}{(s+1)^2} ds  \nonumber
+ \fr{C}{(t+1)^r} \int_{0}^{t} \fr{(s+1)^r}{(s+1)^3} \left(\int_{0}^{s} \|u(\tau)\|_{L^2} ^2 d\tau \right)^2 ds   
\ee for all $t \ge 0$. Here $\Gamma_9$ is a constant depending only on the initial data and $r$.
We have 
\be 
\fr{\Gamma_{9}}{(t+1)^{r}} \int_{0}^{t} \frac{(s+1)^{r}}{(s+1)^2} ds
= \fr{\Gamma_{9}}{(r-1)(t+1)^{r}} \left((t+1)^{r-1} - 1 \right)
\le \fr{\Gamma_{9}}{(r-1)(t+1)} 
\ee for any $r > 1$. 
Moreover, applying the Cauchy-Schwarz inequality in the time variable yields
\be 
\left(\int_{0}^{s} \|u(\tau)\|_{L^2}^2 d\tau \right)^2 \le s \int_{0}^{s} \|u(\tau)\|_{L^2}^4 d\tau, 
\ee so that 
\beg{align} 
\fr{C}{(t+1)^r} \int_{0}^{t} \fr{(s+1)^r}{(s+1)^3} \left(\int_{0}^{s} \|u(\tau)\|_{L^2} ^2 d\tau \right)^2 ds
&\le \fr{C}{(t+1)^r} \left(\int_{0}^{t} (s+1)^{r-2} ds\right) \left(\int_{0}^{t} \|u(s)\|_{L^2}^4 ds \right) \nonumber
\\&\le \fr{C}{(r-1)(t+1)} \left(\int_{0}^{t} \|u(s)\|_{L^2}^4 ds \right)
\end{align} for any $r > 1$. Taking $r =2$ and using \eqref{EQ4} give 
\be 
\|u(t)\|_{L^2}^2 \le \fr{\Gamma_{10}}{t+1} + \fr{\Gamma_{10}}{t+1} \int_{0}^{t} \fr{\|u(s)\|_{L^2}^2}{\ln (e+s)} ds
\ee and so
\be  
(t+1)\|u(t)\|_{L^2}^2 \le \Gamma_{10} + C'\Gamma_{10} \int_{0}^{t} \fr{(s+1)\|u(s)\|_{L^2}^2}{(s+e) \ln (e+s)} ds
\ee for any $t \ge 0$. By Gronwall's inequality, we obtain 
\beg{align} 
&(t+1) \|u(t)\|_{L^2}^2 
\le \Gamma_{10} + C'\Gamma_{10}^2 \int_{0}^{t} \fr{e^{\int_{s}^{t} \fr{1}{(e + \tau ) \ln (e + \tau)} d\tau}}{(e + s) \ln (e+s)} ds \nonumber
\\&=  \Gamma_{10} + C'\Gamma_{10}^2 \int_{0}^{t} \fr{\ln (e+t)}{(e + s) \left[\ln (e+s)\right]^2} ds
\le \Gamma_{10} + C'\Gamma_{10}^2 \ln (e + t).
\end{align}  Therefore, 
\be 
\|u(t)\|_{L^2}^2 \le \fr{\Gamma_{11} \ln (t + e)}{t + 1}
\ee for any $t \ge 0$, where $\Gamma_{11}$ is a constant depending only on the initial data. 

\textit{\textbf{Step 6.3} (Sharp decaying bound for the $L^2$ norm of the velocity).} Finally, we prove \eqref{ltime2}. We take $\rho_1(t) = \sqrt{2} (t+1)^{-\fr{1}{2}}$ as in the previous sub-step, and we obtain the bound  
\be \la{abstractbound1}
\|u(t)\|_{L^2}^2 
\le \fr{\|u_0\|_{L^2}^2}{(t+1)^2}
+ \fr{\Gamma_{12}}{t+1}
+ \fr{C}{(t+1)^2} \int_{0}^{t} \fr{1}{s+1} \left(\int_{0}^{s} \|u(\tau)\|_{L^2} ^2 d\tau \right)^2 ds   
\ee for all $t \ge 0$. We note that 
\be
\int_{0}^{s} \|u(\tau)\|_{L^2} ^2 d\tau \le \Gamma_{13} \int_{0}^{s} \fr{\ln (\tau + e)}{\tau + 1}  d\tau
\le C\Gamma_{13}\int_{0}^{s} \fr{\ln (\tau + e)}{\tau + e}  d\tau 
\le \Gamma_{14} \left[\ln (s + e)\right]^2
\ee and so
\be 
\int_{0}^{t} \fr{1}{s+1} \left(\int_{0}^{s} \|u(\tau)\|_{L^2} ^2 d\tau \right)^2 ds 
\le \Gamma_{15} \int_{0}^{t} \fr{1}{\sqrt{s+1}} ds 
\ee for all $t \ge 0$. Therefore,
\be 
\|u(t)\|_{L^2}^2 \le \fr{\Gamma_{16}}{t+1}
\ee for all $t \ge 0$, where $\Gamma_{16}$ is a positive constant depending only on the initial data. This ends the proof of Theorem \ref{decay}.

Now we study the rate of convergence of the gradients of the charge density and the velocity. 

\beg{thm} \la{decay2} Let $u_0 \in H^1 \cap L^1$ be divergence-free such that $\na u_0 \in L^1$. Let $q_0 \in H^1 \cap L^1$ such that $\na q_0 \in L^1$. There exist positive constants $K_0$ and $K_0'$ depending only on the initial data and some universal constants such that the unique global-in-time solution $(q,u)$ of \eqref{intro1}--\eqref{qinit} obeys 
\be \la{lt1}
\|\na u(t)\|_{L^2}^2 \le \fr{K_0}{t+1}
\ee and
\be \la{lt2}
\|\na q(t)\|_{L^2}^2 \le \fr{K_0'}{(t + 1)^{2}}
\ee 
for all $t \ge 0$. 
\end{thm}

\textbf{Proof:} The proof is divided into 5 steps.

\textit{\textbf{Step 1} (Pointwise bounds for the Fourier transform of $\na u$).} The Fourier transform of the gradient of $u$ satisfies 
\be 
\pa_t \widehat{\na u} + \widehat{\na \PP(u \cdot \na u)} - \widehat{\na \Delta u} = - \widehat{\na \PP(qRq)}
\ee yielding the differential inequality
\be 
\pa_t \widehat{\na u} + |\xi|^2 \widehat{\na u} \le |\xi|^2 \|u\|_{L^2}^2 + |\xi| \|q\|_{L^2}^2.
\ee Integrating in time from $0$ to $t$, and using \eqref{step11}, we obtain
\be \la{decay2st1}
|\widehat{\na u} (\xi,t)| \le \|\na u_0\|_{L^1} + K_1|\xi| + K_2 |\xi|^2 t
\ee where $K_1$ and $K_2$ are positive constants depending only on the initial data. 

We note that the pointwise estimate \eqref{decay2st1} is not the sharpest. Indeed, one can use the decaying bounds for the $L^2$ norms of the velocity and the surface charge density derived in Theorem \ref{decay} instead of using \eqref{step11}, and this will result in a bound whose growth in $t$ is slower when $t \ge 1$. However, this will not improve the decay in the upcoming step, so we disregard this observation. 

\textit{\textbf{Step 2} (Decaying bound for the $L^2$ norm of $\na u$).}
We take the $L^2$ inner product of equation \eqref{intro2} with $-\Delta u$ and we get
\be \la{gradu}
\fr{1}{2} \fr{d}{dt} \|\na u\|_{L^2}^2 + \|\Delta u\|_{L^2}^2 =  \int_{\R^2} qRq \cdot \Delta u.
\ee The nonlinear term $(u \cdot \na u, \Delta u)_{L^2}$ vanishes due to the fact that the matrix $M^tM^2$ has a zero  trace where $M$ is the two-by-two traceless matrix whose entries are given by $M_{ij} = \fr{\pa u_i}{\pa x_j}$ and $M^t$ is its transpose.  
In view of H\"older's inequality with exponents 4,4,2, the boundedness of the Riesz transforms on $L^4(\RR^2)$ and Young's inequality, we obtain 
\be \la{difna}
\fr{d}{dt} \|\na u\|_{L^2}^2 + \|\Delta u\|_{L^2}^2 \le C\|q\|_{L^4}^4.
\ee
Using the $L^4$ estimate \eqref{Bernouilli}, we have 
\be  \la{intl4}
\|q\|_{L^4}^4 \le K_3 (1+t)^{-2}
\ee where $K_3$ depends on the initial data.
We note that the initial charge density is assumed to be in $H^1$ and so it belongs to $L^4$ due to the Sobolev embedding of $H^1(\R^2)$ into $L^4(\R^2)$. 
Going back to \eqref{difna}, we have 
\be 
\fr{d}{dt} \|\na u\|_{L^2}^2 + \|\Delta u\|_{L^2}^2 \le CK_3(t+1)^{-2}.
\ee For $t\in [0, \infty)$, we let 
\be \la{ro2}
\rho_2(t) = r^{\fr{1}{2}}(t+1)^{-\fr{1}{2}}
\ee for some $r>0$ to be chosen later. By Parseval's identity, we get
\be 
\fr{d}{dt} \|\na u\|_{L^2}^2 + \rho_2(t)^2 \|\na u\|_{L^2}^2 
\le CK_3 (1+t)^{-2} + \rho_2(t)^2 \int_{|\xi| \le \rho_2(t)} |\widehat{\na u} (\xi, t)|^2 d\xi. 
\ee
In view of the pointwise bound \eqref{decay2st1}, we have
\beg{align} 
\fr{d}{dt} \|\na u\|_{L^2}^2 + \rho_2(t)^2 \|\na u\|_{L^2}^2 
&\le K_4(t+1)^{-2} + K_5 \rho_2(t)^2 \left[\rho_2(t)^2 + \rho_2(t)^4 + \rho_2(t)^6t^2\right] 
\end{align} for $t \ge 0$.
We multiply by the integrating factor $(t+1)^{r}$ and then we integrate in time from $0$ to $t$. We obtain
\beg{align} 
\|\na u\|_{L^2}^2 
&\le \fr{\|\na u_0\|_{L^2}^2}{(t+1)^{r}} 
+ \fr{K_4}{(t+1)^{r}} \int_{0}^{t} \fr{(s+1)^{r}}{(s+1)^{2}} ds \nonumber
\\&+ \fr{K_6}{(t+1)^{r}} \int_{0}^{t} (1+s)^{r}\left(\fr{1}{(1+s)^2} + \fr{1}{(1+s)^3} \right) ds 
\end{align} where $K_6$ depends only on the initial data. We can take any $r > 2$ and we obtain the bound \eqref{lt1}.

\textit{\textbf{Step 3} (Bounds for $\int_{0}^{t} (s+1)^{\gamma}\|\Delta u(s)\|_{L^2}^2 ds$ where $\gamma\ne 1$ is a real number).} 
Let $\gamma \ne 1$. The differential inequality \eqref{difna} yields 
\be 
\fr{d}{dt} ((t+1)^{\gamma} \|\na u\|_{L^2}^2) - \gamma (t+1)^{\gamma -1} \|\na u\|_{L^2}^2 + (t+1)^{\gamma} \|\Delta u\|_{L^2}^2 \le C(t+1)^{\gamma} \|q\|_{L^4}^4
\ee for all $t \ge 0$. Integrating in time from $0$ to $t$ and using \eqref{Bernouilli} and \eqref{lt1}, we obtain
\beg{align} \la{disb}
\int_{0}^{t} (s+1)^{\gamma} \|\Delta u\|_{L^2}^2 ds 
&\le \|\na u_0\|_{L^2}^2 + (\gamma K_0 + K_7) \int_{0}^{t} (s+1)^{\gamma - 2} ds  \nonumber
\\&= \|\na u_0\|_{L^2}^2 + \fr{\gamma K_0 + K_7}{\gamma - 1} \left[(t+1)^{\gamma -1} - 1 \right] 
\end{align} for some positive constant $K_7$ depending on $\|q_0\|_{L^2}$ and $\|q_0\|_{L^4}$.

\textit{\textbf{Step 4} (Pointwise bounds for the Fourier transform of $\na q$).}
The Fourier transform of the gradient of $q$ satisfies
\be  
\pa_t \widehat{\na q} + \widehat{\na (u \cdot \na q)} + \widehat{\na \l q} = 0,
\ee hence
\be 
\pa_t \widehat{\na q} + |\xi| \widehat{\na q} \le |\xi|^2 \|u\|_{L^2}\|q\|_{L^2}.
\ee
Using \eqref{step11}, we obtain the pointwise bound
\be \la{naqf}
|\widehat{\na q} (\xi, t)| \le \|\na q_0\|_{L^1} + K_8 |\xi|^2 \sqrt{t} 
\ee for all $\xi \in \R^2$ and $t \ge 0$. Here $K_8 $ depends only on the initial data. 

\textit{\textbf{Step 5} (Decaying bound for the $L^2$ norm of $\na q$).}
The $L^2$ norm of the gradient of $q$ evolves according to the energy equality
\be 
\fr{1}{2} \fr{d}{dt} \|\na q\|_{L^2}^2 + \|\l^{\fr{3}{2}}q\|_{L^2}^2 = (u \cdot \na q, \Delta q)_{L^2}.
\ee
In view of the Ladyzhenskaya interpolation inequality 
\be 
\|\na u\|_{L^4} \le C\|\na u\|_{L^2}^{\fr{1}{2}} \|\Delta u\|_{L^2}^{\fr{1}{2}}
\ee and the interpolation inequality \cite{AI}
\be 
\|\na q\|_{L^{\fr{8}{3}}}^2 \le C\|q\|_{L^4}^{\fr{1}{2}} \|\l^{\fr{3}{2}}q\|_{L^2}^{\fr{3}{2}},
\ee we estimate the nonlinear term
\be 
|(u \cdot \na q, \Delta q)_{L^2}| \le \|\na u\|_{L^4} \|\na q\|_{L^{\fr{8}{3}}}^2 
\le C\|\na u\|_{L^2}^{\fr{1}{2}} \|\Delta u\|_{L^2}^{\fr{1}{2}} \|q\|_{L^4}^{\fr{1}{2}} \|\l^{\fr{3}{2}}q\|_{L^2}^{\fr{3}{2}}.
\ee
Applying Young's inequality, we obtain
\be \la{thtwo} 
\fr{d}{dt} \|\na q\|_{L^2}^2 + \|\l^{\fr{3}{2}} q\|_{L^2}^2 \le C\|\na u\|_{L^2}^2\|\Delta u\|_{L^2}^2\|q\|_{L^4}^2.
\ee
In view of \eqref{lt1} and \eqref{intl4}, we have 
\be 
\fr{d}{dt} \|\na q\|_{L^2}^2 + \|\l^{\fr{1}{2}} \na q\|_{L^2}^2 \le \fr{K_{9}}{(t+1)^{2}} \|\Delta u\|_{L^2}^2
\ee for all $t \ge 0$. 	Here $K_{9}$ depends only the initial data. Letting 
\be 
\rho_3(t)= r(t+1)^{-1},
\ee we split the dissipation term,
\be 
\|\l^{\fr{3}{2}} q\|_{L^2}^2  \ge \rho_3(t) \|\na q\|_{L^2}^2 - \rho_3(t) \int_{|\xi| \le \rho_3(t)} |\widehat{\na q} (\xi, t)|^2 d\xi
\ee yielding the differential inequality
\be 
\fr{d}{dt} \|\na q\|_{L^2}^2 + \rho_3(t) \|\na q\|_{L^2}^2 \le \fr{K_{9}}{(t+1)^{2}} \|\Delta u\|_{L^2}^2 + \rho_3(t) \int_{|\xi| \le \rho_3(t)} |\widehat{\na q} (\xi, t)|^2 d\xi
\ee
In view of the pointwise bound for $\widehat{\na q}$ given by \eqref{naqf}, we have 
\be 
\int_{|\xi| \le \rho_3(t)} |\widehat{\na q} (\xi, t)|^2 d\xi
\le K_{10} \left(\rho_3(t)^2 + \rho_3(t)^6 t \right)
\ee and so
\be 
\fr{d}{dt} \|\na q\|_{L^2}^2 + \rho_3(t) \|\na q\|_{L^2}^2 \le \fr{K_{9}}{(t+1)^{2}} \|\Delta u\|_{L^2}^2 + K_{10}\left(\rho_3(t)^3 + \rho_3(t)^7t \right). 
\ee
We multiply both sides by $(t+1)^{r}$ and we integrate in time from $0$ to $t$. We obtain
\begin{align} \la{grq1}
\|\na q(t)\|_{L^2}^2 &\le \fr{\|\na q_0\|_{L^2}^2}{(t+1)^{r }}
+ \fr{K_{11}}{(t+1)^{r}} \int_{0}^{t} (s+1)^{r-2}\|\Delta u(s)\|_{L^2}^2 ds\\
&+ \fr{K_{12}}{(t+1)^{r}} \int_{0}^{t} \left[(s+1)^{r-3} - (s+1)^{r-6}\right] ds.
\end{align}
In view of \eqref{disb} applied with $\gamma = r-2$, we have
\be 
\int_{0}^{t} (s+1)^{r-2}\|\Delta u(s)\|_{L^2}^2 ds
\le  \|\na u_0\|_{L^2}^2 + \fr{(r -2)K_0 + K_7}{r - 3} \left[(t+1)^{r-3} - 1 \right] 
\ee for any $r \ne 3$, and so
\be \la{grq2}
\fr{K_{11}}{(t+1)^{r}} \int_{0}^{t} (s+1)^{r-2}\|\Delta u(s)\|_{L^2}^2 ds
\le \fr{K_{13}}{(t+1)^3}
\ee for any $r > 3$. Here $K_{13}$ depends on the initial data and $r$. Putting \eqref{grq1} and \eqref{grq2} together and choosing $r = 6$ give the desired decay \eqref{lt2}. This completes the proof of Theorem \ref{decay2}.

Now we establish decaying bounds for higher order derivatives. We need the following proposition.

\beg{prop} \la{threetwo} Let $u_0 \in H^1 \cap L^1$ be divergence-free such that $\na u_0 \in L^1$. Let $q_0 \in H^1 \cap L^1$ such that $\na q_0 \in L^1$. Let $\beta > 3$. There exist a positive universal constant $C$ and positive constants $c_1$, $c_2$, and $c_3$ depending only on the initial data such that the solution $q$ of \eqref{intro1}--\eqref{qinit} obeys 
\be  \la{tht2}
\int_{0}^{t} (s+1)^{\beta} \|\l^{\fr{3}{2}}q(s)\|_{L^2}^2 ds 
\le \|\na q_0\|_{L^2}^2 + C\|\na u_0\|_{L^2}^2 + C\fr{(\beta - 2)c_1 + c_2}{\beta - 3} (t+1)^{\beta - 3}
+ \fr{\beta c_3}{\beta - 2} (t+1)^{\beta - 2}
\ee for all $t \ge 0$. 
\end{prop}

\textbf{Proof:} In view of the differential inequality \eqref{thtwo}, we have
\be 
\fr{d}{dt} (t+1)^{\beta} \|\na q\|_{L^2}^2 -\beta(t+1)^{\beta -1} \|\na q\|_{L^2}^2 + (t+1)^{\beta} \|\l^{\fr{3}{2}}q\|_{L^2}^2 \le C(t+1)^{\beta} \|\na u\|_{L^2}^2 \|\Delta u\|_{L^2}^2 \|q\|_{L^4}^2.
\ee Integrating in time from $0$ to $t$, using the bounds \eqref{Bernouilli} and \eqref{lt1} and applying \eqref{disb} with $\gamma = \beta -2$, we obtain \eqref{tht2}.

\beg{thm} \la{decay3} Let $u_0 \in H^2 \cap L^1$ be divergence-free such that $\na u_0 \in L^1$ and $\Delta u_0 \in L^1$. Let $q_0 \in H^2 \cap L^1$ such that $\na q_0 \in L^1$ and $\Delta q_0 \in L^1$. There exist positive constants $M_0$ and $M_0'$ depending only on the initial data and some universal constants such that the unique global-in-time solution $(q,u)$ of \eqref{intro1}--\eqref{qinit} obeys 
\be \la{deltau1}
\|\Delta u(t)\|_{L^2}^2 \le \fr{M_0}{t+1}
\ee and
\be \la{deltaq2}
\|\Delta q(t)\|_{L^2}^2 \le \fr{M_0'}{(t + 1)^2}
\ee 
for all $t \ge 0$. 
\end{thm}

\textbf{Proof:} The Fourier transform of $\Delta u$ obeys
\be 
\pa_t \widehat{\Delta u} + |\xi|^2 \widehat{\Delta u}
\le C|\xi|^3 \|u\|_{L^2}^2 + C|\xi|^2 \|q\|_{L^2}^2.
\ee Consequently, it satisfies the pointwise bound 
\be \la{ptb}
|\widehat{\Delta u} (\xi, t)| \le \|\Delta u_0\|_{L^1} + M_1 |\xi|^3 t + M_2 |\xi|^2
\ee for all $\xi \in \R^2$ and all $t \ge 0$. Here $M_1$ and $M_2$ are positive constants depending only on the initial data. 
The $L^2$ norm of $\Delta u$ evolves according to the energy equality 
\be \la{deltaeq1}
\fr{1}{2} \fr{d}{dt} \|\Delta u\|_{L^2}^2 + \|\na \Delta u\|_{L^2}^2
= - \int_{\RR^2} \Delta (qRq) \cdot \Delta u dx 
- \int_{\RR^2} \Delta (u \cdot \na u) \cdot \Delta u dx.
\ee Integrating by parts, using \eqref{intro3}, and applying Ladyzhenskaya's interpolation inequality, we estimate the second term on the right hand side in \eqref{deltaeq1} as 
\be 
\left|\int_{\RR^2} \Delta (u \cdot \na u) \cdot \Delta u dx \right|
\le C\|\na u\|_{L^2}\|\Delta u\|_{L^2}\|\na \Delta u\|_{L^2}.
\ee In view of the boundedness of the Riesz transforms on $L^4$ and the continuous embedding of $\dot{H}^{\fr{1}{2}}$ in $L^4$, we obtain for the first term on the right hand side in \eqref{deltaeq1} 
\be \la{deltaeq2}
\left|\int_{\RR^2} \Delta (qRq) \cdot \Delta u dx  \right|
\le C\|q\|_{L^4}\|\l^{\fr{3}{2}}q\|_{L^2}\|\na \Delta u\|_{L^2}.
\ee
From \eqref{deltaeq1}--\eqref{deltaeq2} and using Young's inequality, we obtain the energy inequality 
\be 
\fr{d}{dt} \|\Delta u\|_{L^2}^2 + \|\na \Delta u\|_{L^2}^2 
\le C\|q\|_{L^4}^2 \|\l^{\fr{3}{2}}q\|_{L^2}^2 + C\|\na u\|_{L^2}^2 \|\Delta u\|_{L^2}^2.
\ee
In view of Parseval's identity, we have 
\be 
\fr{d}{dt} \|\Delta u\|_{L^2}^2 + \rho_2(t)^2 \|\Delta u\|_{L^2}^2
\le C\|q\|_{L^4}^2 \|\l^{\fr{3}{2}}q\|_{L^2}^2 
+ C\|\na u\|_{L^2}^2 \|\Delta u\|_{L^2}^2 
+ \rho_2(t)^2 \int_{|\xi| \le \rho_2(t)} |\widehat{\Delta u}(\xi, t)|^2 d\xi
\ee
where $\rho_2$ is the function defined by \eqref{ro2}. 
The decay bounds \eqref{intl4} and \eqref{lt1} together with the pointwise bound for the Fourier transform of $\Delta u$ given by \eqref{ptb} yield 
\beg{align} 
\fr{d}{dt} \|\Delta u\|_{L^2}^2 + \rho_2(t)^2 \|\Delta u\|_{L^2}^2 
\le \fr{M_3}{t+1} \|\l^{\fr{3}{2}}q\|_{L^2}^2 
+ \fr{M_4}{t+1} \|\Delta u\|_{L^2}^2 
+ M_5 \rho_2(t)^2 \left[\rho_2(t)^2 + \rho_2(t)^8 t^2 + \rho_2(t)^6 \right]
\end{align} for all $t \ge 0$, where $M_3, M_4$ and $M_5$ are positive constants depending only on the initial data. Multiplying by the integrating factor and integrating in time from $0$ to $t$, we obtain
\beg{align} 
&\|\Delta u\|_{L^2}^2 
\le \fr{\|\Delta u_0\|_{L^2}^2}{(t+1)^{r}}
+ \fr{M_3}{(t+1)^{r}} \int_{0}^{t} (s+1)^{r - 1}\|\l^{\fr{3}{2}}q\|_{L^2}^2 ds \nonumber
\\&+ \fr{M_4}{(t+1)^{r}} \int_{0}^{t} (s+1)^{r-1} \|\Delta u\|_{L^2}^2 ds 
+ \fr{M_5}{(t+1)^{r }} \int_{0}^{t} \left[(s+1)^{r-2} + (s+1)^{r-3} +(s+1)^{r-4}\right] ds.
\end{align}
We choose $r = 5$. In view of the bound \eqref{disb} applied with $\gamma = r-1$ and Proposition \ref{threetwo} applied with $\beta = r-1$, we obtain \eqref{deltau1}. 

The Fourier transform of the Laplacian of $q$ satisfies
\be 
\pa_t \widehat{\Delta q} + |\xi| \widehat{\Delta q} \le |\xi|^3 \|u\|_{L^2}\|q\|_{L^2}
\ee and hence
\be \la{Deltaqf}
|\widehat{\Delta q} (\xi, t)| \le \|\Delta q_0\|_{L^1} + M_8 |\xi|^3 \sqrt{t} 
\ee for all $\xi \in \R^2$ and $t \ge 0$. Here $M_8 $ depends only on the initial data. 
Now, we establish decaying estimate for $\|\Delta q\|_{L^2}^2$ which evolves according to 
\be 
\fr{1}{2} \fr{d}{dt} \|\Delta q\|_{L^2}^2 + \|\l^{\fr{5}{2}}q\|_{L^2}^2
= 2\int_{\RR^2} (\na u \cdot \na (\na q)) \Delta q dx
+ \int_{\RR^2} (\Delta u \cdot \na q) \Delta q dx.
\ee In view of the Gagliardo-Nirenberg interpolation inequality
\be 
\|\Delta q\|_{L^2} \le C\|\l^{\fr{5}{2}} q\|_{L^2}^{\fr{4}{5}}\|q\|_{L^2}^{\fr{1}{5}},
\ee the Sobolev embedding inequality
\be 
\|\Delta q\|_{L^4} \le C\|\l^{\fr{5}{2}}q\|_{L^2},
\ee and the bound
\be 
\|\na \na q\|_{L^2} = \|\na \l^{-1} \na \l^{-1} \Delta q\|_{L^2} \le C\|\Delta q\|_{L^2}
\ee that follows from the boundedness of the Riesz transforms on $L^2$, we obtain 
\beg{align} 
\fr{1}{2} \fr{d}{dt} \|\Delta q\|_{L^2}^2 + \|\l^{\fr{5}{2}}q\|_{L^2}^2 
&\le C\|\na u\|_{L^4} \|\Delta q\|_{L^2} \|\Delta q\|_{L^4} + C\|\Delta u\|_{L^2}\|\na q\|_{L^4} \|\Delta q\|_{L^4} \nonumber
\\&\le C\|\na u\|_{L^2}^{\fr{1}{2}}\|\Delta u\|_{L^2}^{\fr{1}{2}} \|\l^{\fr{5}{2}} q\|_{L^2}^{\fr{9}{5}} \|q\|_{L^2}^{\fr{1}{5}} + C\|\Delta u\|_{L^2} \|\l^{\fr{3}{2}}q\|_{L^2}\|\Delta q\|_{L^4} \nonumber
\\&\le \fr{1}{2} \|\l^{\fr{5}{2}}q\|_{L^2}^2 
+ C\|\na u\|_{L^2}^5 \|\Delta u\|_{L^2}^5\|q\|_{L^2}^2 
+ C\|\Delta u\|_{L^2}^{2} \|\l^{\fr{3}{2}}q \|_{L^2}^2.
\end{align} Consequently, 
\beg{align} 
&\fr{d}{dt} \|\Delta q\|_{L^2}^2 + \rho_2(t)^2 \|\Delta q\|_{L^2}^2 
\le  C\|\na u\|_{L^2}^5 \|\Delta u\|_{L^2}^5\|q\|_{L^2}^2 
\\&\quad\quad+ C\|\Delta u\|_{L^2}^{2} \|\l^{\fr{3}{2}}q \|_{L^2}^2 \nonumber
+ \rho_2(t)^2 \int_{|\xi| \le \rho_2(t)^2} |\widehat{\Delta q} (\xi,t)|^2 d\xi
\end{align} where $\rho_2$ is defined by \eqref{ro2}.
In view of the estimates \eqref{ltime1}, \eqref{lt1} and \eqref{deltau1}, Proposition \ref{threetwo} applied with $\beta = r-1$, and the pointwise bound for the Fourier transform of $\Delta q$ given by \eqref{Deltaqf}, we obtain \eqref{deltaq2}. This ends the proof of Theorem \ref{decay3}.
 
Let $C^{0,\fr{1}{2}}$ be the space of bounded $1/2$-H\"older continuous functions on $\mathbb{R}^2$ with 
\be 
\|f\|_{C^{0,\fr{1}{2}}} = \|f\|_{L^{\infty}} + \sup\limits_{x, y \in \mathbb{R}^2, x 
\neq y} \frac{|f(x) - f(y)|}{|x-y|^{\fr{1}{2}}}.
\ee
In view of the continuous Sobolev embedding of $W^{1,4}$ into $C^{0, \fr{1}{2}}$, the Ladyzhenskaya interpolation inequality, and Theorems \ref{decay}, \ref{decay2}, and \ref{decay3}, we obtain the following statement.

\beg{cor}
Let $u_0 \in H^2 \cap L^1$ be divergence-free such that $\na u_0 \in L^1$ and $\Delta u_0 \in L^1$. Let $q_0 \in H^2 \cap L^1$ such that $\na q_0 \in L^1$ and $\Delta q_0 \in L^1$. There exist positive constants $A_0$ and $A_0'$ depending only on the initial data and some universal constants such that the unique global-in-time solution $(q,u)$ of \eqref{intro1}--\eqref{qinit} obeys 
\be \la{hold1}
\|u(t)\|_{C^{0,\fr{1}{2}}}^2 \le \fr{A_0}{t+1}
\ee and
\be \la{hold2}
\|q(t)\|_{C^{0, \fr{1}{2}}}^2 \le \fr{A_0'}{(t + 1)^{2}}
\ee 
for all $t \ge 0$. 
\end{cor}

\section{Decomposition of the solution} \la{s3}

In this section, we decompose the charge density $q$ and the velocity $u$ solutions of \eqref{intro1}--\eqref{qinit} in the sum of solutions $Q$ and $U$ of the linear equations 
\be \la{linear}
\pa_t Q + \l Q = 0
\ee and
\be \la{linear2}
\pa_t U - \Delta U = 0
\ee with initial datum $Q(0) = q_0$ and $U(0) =u_0$ and remainders. We study the decays of the remainders  $q-Q$ and $u - U$ in $L^2$ and we show that they are faster than the decays of the $L^2$ norms of $q$ and $u$ respectively.
The solutions of \eqref{linear} and \eqref{linear2} are given explicitly by 
\be 
Q(t) = \int_{\R^2} K_t^1(x-w) q_0(w) dw 
\ee 
and 
\be 
U(t) = \int_{\R^2} K_t^2(x-w) u_0(w) dw
\ee
where $K_t^s$ is the kernel defined by its Fourier transform 
\be 
\mathcal{F} (K_t^s)(\xi) = e^{-|\xi|^st}.
\ee

The following proposition describes the decay of $\na Q$ and $\na U$ in $L^{\infty}$.

\beg{prop} Suppose $q_0 \in L^1$ such that $\int_{\RR^2} |\xi| |\widehat{q_0}(\xi)| d\xi < \infty $ and $u_0 \in L^1$ such that $\int_{\RR^2} |\xi| |\widehat{u_0}(\xi)| d\xi < \infty $. Then there exist positive constants $R_0$ and $R_0'$ depending only on the initial data such that the solutions $Q$ and $U$ of the linear equations \eqref{linear} and \eqref{linear2} satisfy
\be \la{linearinf}
\|\na Q(t)\|_{L^{\infty}} \le \frac{R_0}{(t+1)^3}
\ee and 
\be \la{linearinf2}
\|\na U(t)\|_{L^{\infty}} \le \fr{R_0'}{(t+1)^{\fr{3}{2}}}
\ee for all $t \ge 0$. 
\end{prop}

\textbf{Proof:} In view of Parseval's identity and the translation property of the Fourier transform, we have
\be 
\na Q (x) = \na K_t^1 * q_0 (x) = \int_{\R^2} \na K_t^1 (x-w) q_0(w) dw = \int_{\RR^2} e^{-2\pi i x \cdot \xi} \widehat{\na K_t^1}(\xi) \widehat{q_0}(\xi) d\xi.
\ee On one hand,
\be 
\|\na Q\|_{L^{\infty}} \le C\|\widehat{q_0}\|_{L^{\infty}} \int_{\RR^2} |\xi| |\widehat{K_t^1}(\xi)| d\xi
\le C\|q_0\|_{L^1} \int_{0}^{\infty} r^2 e^{-rt} dr \le C\|q_0\|_{L^1} t^{-3}
\ee and so
\be 
t^3 \|\na Q\|_{L^{\infty}} \le C\|q_0\|_{L^1}. 
\ee On the other hand, 
\be 
\|\na Q\|_{L^{\infty}} \le C \int_{\RR^2} |\xi| |\widehat{K_t^1}(\xi)||\widehat{q_0}(\xi)| d\xi
\le C \int_{\RR^2} |\xi| |\widehat{q_0}(\xi)| d\xi
\ee for all $t \ge 0$. Hence 
\be 
(1+t)^3 \|\na Q\|_{L^{\infty}} \le 4(1+t^3)\|\na Q\|_{L^{\infty}} \le C\left(\|q_0\|_{L^1} + \int_{\RR^2} |\xi| |\widehat{q_0}(\xi)| d\xi \right) 
\ee for all $t \ge 0$, yielding \eqref{linearinf}. Similarly, we have 
\be 
(1+t)^{\fr{3}{2}} \|\na U\|_{L^{\infty}} \le C(1+t^{\fr{3}{2}})\|\na U\|_{L^{\infty}} \le C\left(\|u_0\|_{L^1} + \int_{\RR^2} |\xi| |\widehat{u_0}(\xi)| d\xi \right) 
\ee for all $t \ge 0$, yielding \eqref{linearinf2}.

\beg{rem} The assumptions $\int_{\RR^2} |\xi| |\widehat{q_0}(\xi)| d\xi < \infty $ and $\int_{\RR^2} |\xi| |\widehat{u_0}(\xi)| d\xi < \infty $ are required to obtain the uniform-in-time boundedness of the $L^{\infty}$ norms $\na Q$ and $\na U$ for small times $t \in (0,1)$. This imposed regularity can be dropped since we are interested in studying the long-time behavior of solutions.  
\end{rem}

Next, we consider the pointwise behavior of the Fourier transforms of the differences $q-Q$ and $u- U$. We need first the following lemmas.

\beg{lem} \la{lem1} For $f \in L^2(\RR^2)$ and $x \in \RR^2$, we let 
\be 
Tf(x) = \lim\limits_{\epsilon \to 0} \int_{|x-y| > \epsilon} \fr{\sqrt{|y|^2 +1} - \sqrt{|x|^2 +1}}{|x-y|^{3}} f(y)dy.
\ee There exists a universal constant $C>0$ (independent of $f$) such that 
\be 
\|Tf\|_{L^2} \le C\|f\|_{L^2}.
\ee
\end{lem}

\textbf{Proof:} We write
\be 
Tf(x) = \lim\limits_{\epsilon \to 0} \int_{|x-y| > \epsilon} (a(y) - a(x))k(x-y) f(y)dy.
\ee where $a(x)$ is the function defined on $\RR^2$ by
\be 
a(x) = \sqrt{|x|^2 + 1}
\ee and $k(x)$ is the function defined on $\RR^2 \setminus \left\{0\right\}$ by 
\be 
k(x) = \fr{1}{|x|^3}.
\ee We note that $k$ is homogeneous of degree $-3$. Moreover, the gradient of $a$ is given by 
\be 
\na a(x) = \left(\fr{x_1}{\sqrt{|x|^2 +1}}, \fr{x_2}{\sqrt{|x|^2 + 1}} \right)
\ee and satisfies $\|\na a\|_{L^{\infty}} \le 1$. Therefore, $T$ is a well-defined operator and bounded on $L^2$ (see page 435 in Section 2 of \cite{C}).

Using Lemma \ref{lem1}, we study the evolution of $(\sqrt{|x|^2 + 1})q(x)$ in $L^2(\R^2)$.

\beg{lem} \la{lem2} Let $u_0 \in H^1 \cap L^1$ be divergence-free such that $\na u_0 \in L^1$. Let $q_0 \in H^1 \cap L^1$ such that $\na q_0 \in L^1$. Furthermore, suppose that $\int_{\RR^2} |x|^2q_0(x)^2 dx < \infty$. Then there exists a positive constant $R_1 > 0$ depending only on the initial data such that 
\be \la{moment}
\|(\sqrt{|\cdot|^2 + 1}) q(\cdot,t)\|_{L^2} \le R_1 \ln(t+1) + \|(\sqrt{|\cdot|^2 + 1}) q_0 (\cdot)\|_{L^2} 
\ee holds for all $t \ge 0$.
\end{lem}

\textbf{Proof:} Let $a(x) = \sqrt{|x|^2 + 1}$. The evolution of $aq$ is described by 
\be 
\pa_t (aq) + au \cdot \na q + a\l q = 0.
\ee Multiplying by $aq$ and integrating in the space variable over $\RR^2$, we obtain
\be 
\fr{1}{2} \fr{d}{dt} \|aq\|_{L^2}^2 + \int_{\RR^2} (a \l q) a q 
= - \int_{\RR^2} (a u \cdot \na q) a q .
\ee The cancellation 
\be 
\int_{\RR^2} (u \cdot\na (aq)) aq  = 0 
\ee holds due to \eqref{intro3}, so we can rewrite the nonlinear term as 
\be 
-\int_{\RR^2} (a u \cdot \na q) a q
= \int_{\RR^2} (u \cdot \na a) q^2 a.
\ee By H\"older's inequality, Ladyzhenskaya's interpolation inequality, and the decaying bounds for the $L^2$ norms of $q, u, \na u$ and $\na q$ given by \eqref{ltime1}, \eqref{ltime2}, \eqref{lt1} and \eqref{lt2}, respectively,  we estimate 
\beg{align} 
&\left|\int_{\RR^2} (u \cdot \na a) q^2 a \right|
\le \|\na a\|_{L^{\infty}} \|q\|_{L^4} \|u\|_{L^4} \|aq\|_{L^2} \nonumber
\\&\quad\quad\le C\|q\|_{L^2}^{\fr{1}{2}} \|\na q\|_{L^2}^{\fr{1}{2}} \|u\|_{L^2}^{\fr{1}{2}}\|\na u\|_{L^2}^{\fr{1}{2}}\|aq\|_{L^2}
\le R_2 (t+1)^{-\fr{3}{2}} \|aq\|_{L^2}
\end{align} for some constant $R_2$ depending only on the initial data. 
Now we write the linear term as the sum
\beg{align} 
\int_{\RR^2} (a \l q) a q
&= \int_{\RR^2} aq \l (aq)
+ \int_{\RR^2} (aq) \left[a\l q - \l (aq) \right] \nonumber
\\&= \|\l^{\fr{1}{2}} (aq)\|_{L^2}^2 + \int_{\RR^2} (aq) \left[a\l q - \l (aq) \right] .
\end{align} By the Cauchy-Schwarz inequality, we bound
\be 
\left|\int_{\RR^2} (aq) \left[a\l q - \l (aq) \right]\right| \le \|aq\|_{L^2} \|a\l q - \l(aq)\|_{L^2}.
\ee The pointwise formula for the fractional Laplacian of order 1 yields
\beg{align} 
(a \l q - \l (aq))(x) 
&= C\int_{\RR^2} \left[\fr{a(x)q(x) - a(x) q(y)}{|x-y|^3} - \fr{a(x) q(x) - a(y) q(y)}{|x-y|^3} \right] dy \nonumber
\\&= C\int_{\RR^2} \fr{a(y) - a(x)}{|x-y|^3} q(y) dy
\end{align} where $C$ is  positive universal constant. As a consequence of Lemma \ref{lem1} and \eqref{ltime1}, we obtain
\be 
\|a\l q - \l(aq)\|_{L^2} \le C\|q\|_{L^2}\le  C(t+1)^{-1}.
\ee
Therefore, the $L^2$ norm of $aq$ obeys the energy inequality
\be 
\fr{1}{2} \fr{d}{dt} \|aq\|_{L^2}^2 + \|\l^{\fr{1}{2}} (aq)\|_{L^2} \le \left[R_2(t+1)^{-\fr{3}{2}} + C(t+1)^{-1}\right] \|aq\|_{L^2}
\ee so
\be 
\fr{1}{2} \fr{d}{dt} \|aq\|_{L^2}^2 \le R_3(t+1)^{-1} \|aq\|_{L^2}
\ee for some positive constant $R_3$ depending only on the initial data. Dividing both sides of the inequality by $\|aq\|_{L^2}$, we get 
\be 
\fr{d}{dt} \|aq\|_{L^2} \le R_3 (t+1)^{-1}.
\ee Integrating in time from $0$ to $t$, we obtain \eqref{moment}.

The following lemma is needed to obtain a growth in $|\xi|$ for the Fourier transform of $\PP(qRq)$.

\beg{lem} \la{lem3} Let $f \in L^2(\RR^2)$ such that $\int_{\RR^2} |x|^2f(x)^2 dx < \infty$. Then 
\be \la{lemma}
|\widehat{\PP(fRf)} (\xi)| \le C|\xi| \|f\|_{L^2} \left(\int_{\RR^2} |x|^2 |f(x)|^2 dx\right)^{\fr{1}{2}}.
\ee where $\PP$ is the Leray projector and $R = (R_1, R_2)$ is the Riesz transform vector on $\R^2$.
\end{lem}

\textbf{Proof:} The Leray projector is a Fourier multiplier with a symbol denoted by $m(\xi)$. We have  
\be 
\widehat{\PP(fRf)}(\xi) = m (\xi) \widehat{fRf} (\xi)
\ee for all $\xi \in \RR^2$.
We note that $m(\xi)$ is bounded uniformly in $\xi$. Now, the Fourier transform of $fRf$ at $\xi$ is given by 
\be 
\widehat{fRf} (\xi) = \int_{\RR^2} f(x)Rf(x) e^{- i \xi \cdot x} dx
\ee for $\xi \in \RR^2$. Since the Riesz transform is antisymmetric, we have 
\be 
\int_{\RR^2} f(x)Rf(x) dx = 0
\ee and so we can write $\widehat{fRf}$ at $\xi$ as 
\be  
\widehat{fRf}(\xi)= \int_{\RR^2} f(x)Rf(x) \left(e^{-i \xi \cdot x} -1 \right) dx.
\ee Using the identity
\be 
|e^{-i \xi \cdot x} - 1| \le  |\xi| |x|
\ee that holds for all $x, \xi \in \RR^2$, we estimate
\be 
|\widehat{fRf}(\xi)| \le  |\xi| \int_{\RR^2} |x| |f(x)| |Rf(x)| dx
\le |\xi| \|Rf\|_{L^2} \left(\int_{\RR^2} |x|^2 |f(x)|^2 dx\right)^{\fr{1}{2}}
\ee in view of the Cauchy-Schwarz inequality. This gives the pointwise estimate \eqref{lemma}. 

As a consequence of lemmas \ref{lem2} and \ref{lem3}, we obtain the following statement.

\beg{prop} Let $u_0 \in H^1 \cap L^1$ be divergence-free such that $\na u_0 \in L^1$ and $q_0 \in H^1 \cap L^1$ such that $\na q_0 \in L^1$. Furthermore, suppose that $\int_{\RR^2} |x|^2q_0(x)^2 dx < \infty$. Let $(q,u)$ be the solution of \eqref{intro1}--\eqref{qinit}. Let $\zeta = q - Q$ and $v = u - U$. Then there exist positive constants $R_4$, $R_5$ and $R_6$ depending only on the initial data such that the Fourier transforms of $\zeta$ and $v$ satisfy the pointwise bounds
\be \la{fz}
|\widehat{\zeta} (\xi,t)| \le R_4 |\xi|
\ee and 
\be \la{fv}
|\widehat{v}(\xi, t)| \le R_5 |\xi| \ln (t+1) + R_6 |\xi|\ln^2(t+1)
\ee
for all $\xi \in \R^2$ and $t \ge 0$.
\end{prop}

\textbf{Proof:} The Fourier transform of $\zeta$ obeys
\be 
\pa_t \widehat{\zeta} + |\xi| \widehat{\zeta} = - \widehat{u \cdot \na q} \le |\xi| \|u\|_{L^2}\|q\|_{L^2}.
\ee Consequently, 
\be 
|\widehat{\zeta}(\xi, t)| \le \int_{0}^{t} |\xi| \|u\|_{L^2} \|q\|_{L^2} \le R_4 |\xi|
\ee in view of the decaying bounds \eqref{ltime1} and \eqref{ltime2}.
The Fourier transform of $v$ evolves according to 
\be 
\pa_t \widehat{v} + |\xi|^2 \widehat{v} = - \widehat{\PP(u \cdot \na u)} - \widehat{\PP(qRq)}.
\ee Thus
\be 
|\widehat{v}(\xi,t)| \le C|\xi| \int_{0}^{t} \|u\|_{L^2}^2 ds + C |\xi|\int_{0}^{t} \|q\|_{L^2} \left(\int_{\RR^2}
|x|^2q(x)^2 dx \right)^{\fr{1}{2}} ds
\ee by Lemma \ref{lem3}. In view of Lemma \ref{lem2} and the decaying estimates \eqref{ltime1} and \eqref{ltime2}, we obtain \eqref{fv}.

\beg{thm}\la{nonlin} Let $u_0 \in H^1 \cap L^1$ be divergence-free such that $\na u_0 \in L^1$. Let $q_0 \in H^1 \cap L^1$ such that $\na q_0 \in L^1$. Furthermore, suppose that $\int_{\RR^2} |x|^2q_0(x)^2 dx < \infty$, $\int_{\RR^2} |\xi| |\widehat{q_0}(\xi)| d\xi < \infty $, and $\int_{\RR^2} |\xi| |\widehat{u_0}(\xi)| d\xi < \infty $. Let $(q, u)$ be the solution of \eqref{intro1}--\eqref{qinit}. Then there exist positive constants $R_7$ and $R_8$ depending only on the initial data such that the differences $q - Q$ and $u-U$ satisfy 
\be \la{improved}
\|q(t) - Q(t)\|_{L^2}^2 \le \fr{R_7}{(t+1)^{2 + \fr{3}{2}}}
\ee and 
\be \la{improved2}
\|u(t) - U(t)\|_{L^2}^2 \le \fr{R_8}{(t+1)^{1+ \fr{1}{2}}}
\ee
for all $t \ge 0$.
\end{thm}

\textbf{Proof:} Let $\zeta = q - Q$ and $v = u - U$. We have
\be \la{EQ6}
\pa_t \zeta + \l \zeta = - u \cdot \na q.
\ee Taking the $L^2$ inner product of equation \eqref{EQ6} with $\zeta$ and using \eqref{intro3}, we obtain 
\be 
\fr{1}{2} \fr{d}{dt} \|\zeta\|_{L^2}^2 + \|\l^{\fr{1}{2}} \zeta\|_{L^2}^2 
= \int_{\RR^2} (u \cdot \na q)Q dx
\le \|u\|_{L^2}\|q\|_{L^2}\|\na Q\|_{L^{\infty}}.
\ee As a consequence of Theorem \ref{decay} and the bound \eqref{linearinf}, we obtain the energy inequality
\be 
\fr{d}{dt} \|\zeta\|_{L^2}^2 + \|\l^{\fr{1}{2}} \zeta\|_{L^2}^2
\le \fr{R_9}{(t+1)^{4 + \fr{1}{2}}}
\ee where $R_9$ is a positive constant depending only on the initial data. 
For a fixed $r$, we let $\rho(t) = r(t+1)^{-1}$. Then
\be 
\fr{d}{dt} \|\zeta\|_{L^2}^2 + \rho(t) \|\zeta\|_{L^2}^2
\le \fr{R_9}{(t+1)^{4 + \fr{1}{2}}} 
+ \rho(t) \int_{|\xi| \le \rho(t)} |\widehat{\zeta}(\xi, t)|^2 d\xi.
\ee
Using \eqref{fz}, we estimate
\be
\int_{|\xi| \le \rho(t)} |\widehat{\zeta}(\xi, t)|^2 d\xi
\le R_{10} \rho(t)^4
\ee and we obtain 
\be 
\fr{d}{dt} \|\zeta\|_{L^2}^2 + \rho(t) \|\zeta\|_{L^2}^2
\le  \fr{R_{9}}{(t+1)^{4 + \fr{1}{2} }}  + R_{10} \rho(t)^5.
\ee 
Multiplying by the factor $(s+1)^r$, integrating in the time variable $s$ from $0$ to $t$, and choosing any $r > 4$, we obtain the desired bound \eqref{improved}.  
Now, $v$ obeys 
\be 
\pa_t v - \Delta v = - u \cdot \na u - qRq - \na p.
\ee Taking the $L^2$ inner product of this latter equation with $v$ and using the fact that $v$ is divergence-free, we  get the energy equation 
\be 
\fr{1}{2} \fr{d}{dt} \|v\|_{L^2} + \|\na v\|_{L^2}^2 = \int_{\RR^2} (u \cdot \na u) \cdot U dx
- \int_{\RR^2} (qRq) \cdot v dx.
\ee We estimate
\be 
\int_{\RR^2} (u \cdot \na u) \cdot U dx \le \|u\|_{L^2}^2 \|\na U\|_{L^{\infty}} \le \fr{R_{11}}{(t+1)^{1+ \fr{3}{2}}}
\ee in view of the bounds \eqref{ltime2} and \eqref{linearinf2}, and 
\be 
\int_{\RR^2} (qRq) \cdot v dx \le C\|q\|_{L^4}^2\|v\|_{L^{2}} 
\le C\|q\|_{L^2} \|\na q\|_{L^2} \|v\|_{L^2}
\le \fr{R_{12}}{(t+1)^{1+\fr{3}{2}}}
\ee in view of the decaying estimate \eqref{ltime1}, \eqref{lt2} and \eqref{hold1}. This yields the energy inequality
\be 
\fr{1}{2} \fr{d}{dt} \|v\|_{L^2}^2 + \rho(t) \|v\|_{L^2}^2
\le \fr{R_{13}}{(t+1)^{1+\fr{3}{2}}} + \rho(t) \int_{|\xi| \le \sqrt{\rho(t)}} |\widehat{v}(\xi,t)|^2 d\xi
\ee where $\rho(t) = r(t+1)^{-1}$. Using the pointwise bound for the Fourier transform of $v$ given by \eqref{fv}, we have
\be 
 \int_{|\xi| \le \sqrt{\rho(t)}} |\widehat{v}(\xi,t)|^2 d\xi 
\le R_{14}  \left[\ln^2(t+1) + \ln^4(t+1) \right] \rho(t)^2
\le R_{15} \sqrt{t+1} \rho(t)^2, 
\ee hence
\be 
\fr{1}{2} \fr{d}{dt} \|v\|_{L^2}^2 + \rho(t) \|v\|_{L^2}^2
\le \fr{R_{13}}{(t+1)^{1+\fr{3}{2}}} + R_{15} \sqrt{t+1} \rho(t)^3.
\ee We multiply both sides by $(s+1)^r$, we integrate from $0$ to $t$, we choose any $r > 3/2$, and we obtain \eqref{improved2}.

\section{Appendix: Existence and Uniqueness of Solutions}

In this appendix, we prove the existence of weak and strong solutions for the electroconvection model \eqref{intro1}--\eqref{qinit}. 

\beg{defi}
A solution $(q,u)$ of \eqref{intro1}--\eqref{qinit} is said to be a weak solution on $[0,T]$ if it solves \eqref{intro1}--\eqref{qinit} in the sense of distributions, $u$ is divergence-free in the sense of distributions, 
\be
u \in L^{\infty}(0,T;L^2) \cap L^2(0,T; H^1)
\ee and 
\be
q \in L^{\infty}(0,T; L^2) \cap L^2(0,T; H^{1/2}).
\ee
\end{defi}

\beg{thm} \la{weak} Let $u_0\in L^2$ be divergence-free, let $q_0\in L^2$. Let $T > 0$ be arbitrary. There exists a weak solution $(q,u)$ of the system \eqref{intro1}--\eqref{qinit} on $[0,T]$.
\end{thm}

\textbf{Proof.} We briefly sketch the main ideas of the proof. For $0 < \epsilon \le 1$, we consider a viscous approximation of \eqref{intro1}--\eqref{qinit} given by 
\be \label{sysv}
\begin{cases}  \partial_t q^{\epsilon} + J_{\epsilon}(u^{\epsilon}  \cdot \nabla q^{\epsilon}) + \Lambda q^{\epsilon} - \epsilon \Delta q^{\epsilon} = 0 
\\ \partial_t u^{\epsilon}  + J_{\epsilon}(u^{\epsilon}  \cdot \nabla u^{\epsilon})   - \Delta u^{\epsilon} + \nabla p^{\epsilon}  = -J_{\epsilon} (q^{\epsilon} Rq^{\epsilon})  
\\\nabla \cdot u^{\epsilon} = 0 
\end{cases} 
\ee with smoothed out initial data, where $J_{\epsilon}$ is a standard mollifier operator, $u^{\epsilon} = J_{\epsilon}u$, $q^{\epsilon}= J_{\epsilon}q$ and $p^{\epsilon} = J_{\epsilon}p$.   
For each $\epsilon > 0$, we consider the map
\be 
(q(t), u(t)) \mapsto \Phi_{\epsilon}((q,u))(t) = (e^{\epsilon t \Delta} J_{\epsilon}q_0 - \mathcal{A}^{\epsilon}_t (q^{\epsilon}, u^{\epsilon}), e^{t\Delta}J_{\epsilon}u_0  - \mathcal{B}^{\epsilon}_t (q^{\epsilon}, u^{\epsilon}))
\ee where
\be 
\mathcal{A}^{\epsilon}_t (q^{\epsilon}, u^{\epsilon}) = \int_{0}^{t} e^{\epsilon(t-s)\Delta} J_{\epsilon}(u^{\epsilon} \cdot \na q^{\epsilon})(s) ds + \int_{0}^{t} e^{\epsilon(t-s)\Delta} \l q^{\epsilon}(s) ds
\ee and 
\be 
\mathcal{B}^{\epsilon}_t (q^{\epsilon}, u^{\epsilon}) = \int_{0}^{t} e^{(t-s)\Delta} J_{\epsilon}\PP(u^{\epsilon} \cdot \na u^{\epsilon})(s) ds + \int_{0}^{t} e^{(t-s)\Delta}J_{\epsilon} \PP(q^{\epsilon}Rq^{\epsilon})(s) ds.
\ee
There exists a time $T_{\epsilon} = T_{\epsilon} (\epsilon, \|u_0\|_{L^2}, \|q_0\|_{L^2}) > 0$ such that the map $\Phi_{\epsilon}$ is a contraction on the Banach space 
\be 
X_T = L^{\infty}(0,T; \bar{B}_{L^2} (2\|q_0\|_{L^2}) \oplus L^{\infty}(0,T; \bar{B}_{L^2_{\sigma}} (2\|u_0\|_{L^2}) 
\ee
where $\bar{B}_{L^2}(r)$ is the closed ball in $L^2$, and $\bar{B}_{L^2_{\sigma}}$ is the closed ball in the space of $L^2$ divergence-free vectors. Consequently, $\Phi_{\epsilon}$ has a fixed point $(q^{\epsilon}, u^{\epsilon}) \in X_{T_{\epsilon}}$ solving \eqref{sysv}. This solution extends to the time interval $[0,T]$, and this can be obtained by  establishing uniform-in-time bounds for $(q^{\epsilon}, u^{\epsilon})$ on $[0,T]$. 
Indeed, we have 
\be 
\fr{1}{2} \fr{d}{dt} \left(\|\l^{-\fr{1}{2}}q^{\epsilon} \|_{L^2}^2 + \|u^{\epsilon} \|_{L^2}^2 \right) + \|q^{\epsilon}\|_{L^2}^2 + \|\na u^{\epsilon}\|_{L^2}^2 + \epsilon \|\l^{\fr{1}{2}}q^{\epsilon}\|_{L^2}^2 = 0
\ee as shown in \eqref{weaknormb}. Hence the family of mollified velocities $(u^{\epsilon})_{\epsilon}$ is uniformly bounded in $L^{\infty}(0,T; L^2) \cap L^2(0,T; H^1)$. On the other hand, the $L^2$ norm of $q^{\epsilon}$ evolves according to 
\be 
\fr{1}{2} \fr{d}{dt} \|q^{\epsilon}\|_{L^2}^2 + \|\l^{\fr{1}{2}} q^{\epsilon}\|_{L^2}^2 
+ \epsilon \|\na q^{\epsilon}\|_{L^2}^2 
= 0, 
\ee and so the family of mollified charge densities $(q^{\epsilon})_{\epsilon}$ is uniformly bounded in $L^{\infty}(0,T; L^2) \cap L^2(0,T; H^{\fr{1}{2}})$. The $q^{\epsilon}$ and $u^{\epsilon}$ equations imply that the sequence of time derivatives $(\pa_t q^{\epsilon})_{\epsilon}$ and $(\pa_t u^{\epsilon})_{\epsilon}$ are uniformly bounded in $L^{2}(0,T; H^{-\fr{3}{2}})$ and $L^2(0,T; H^{-1})$ respectively. By the Aubin-Lions lemma, the sequence $((q^{\epsilon}, u^{\epsilon}))_{\epsilon}$ has a subsequence that converges strongly in $L^2(0,T; L^2)$ to a weak solution $(q,u)$ of \eqref{intro1}--\eqref{qinit}. We omit further details.  

\beg{defi} A weak solution $(q,u)$ of \eqref{intro1}--\eqref{qinit} is said to be a strong solution on $[0,T]$ if 
\be
u \in L^{\infty}(0,T;H^1) \cap L^2(0,T; H^2)
\ee
and 
\be 
q \in L^{\infty}(0,T; L^4) \cap L^2(0,T; H^{1/2}).
\ee
\end{defi}

\beg{thm} \la{strong} Let $u_0\in H^1$ be divergence-free and $q_0\in L^4$.  Let $T>0$ be arbitrary. There exists a unique strong solution $(u,q)$ of the system \eqref{intro1}--\eqref{qinit} on $[0,T]$.
\end{thm}

\textbf{Proof.} We take the $L^2$ inner product of the equation satisfied by $q^{\epsilon}$ in   \eqref{sysv} with $(q^{\epsilon})^3$. In view of the divergence-free condition satisfied by $u^{\epsilon}$, the nonlinear term vanishes, that is
\be 
\int_{\RR^2} u^{\epsilon} \cdot \na q^{\epsilon} (q^{\epsilon})^3 dx = 0. 
\ee By the C\'ordoba-C\'ordoba inequality (\cite{CC}), we have 
\be 
\int_{\RR^2} (q^{\epsilon})^3 \l q^{\epsilon} dx  \ge 0 
\ee and 
\be 
-\int_{\RR^2} (q^{\epsilon})^3 \Delta q^{\epsilon} dx \ge 0.
\ee
Consequently, we obtain 
\be 
\fr{1}{4} \fr{d}{dt} \|q^{\epsilon}\|_{L^4}^4  \le 0
\ee which yields the boundedness of $q$ in $L^{\infty} (0, T; L^4(\RR^2))$ by the Banach Alaoglu theorem and the lower semi-continuity of the norm.  
The $L^2$ norm of $\na u^{\epsilon}$ obeys the energy inequality 
\be 
\fr{d}{dt} \|\na u^{\epsilon}\|_{L^2}^2 + \|\Delta u^{\epsilon}\|_{L^2}^2 \le C \|q^{\epsilon}\|_{L^4}^4 
\ee as shown in \eqref{gradu}, yielding the boundedness of $u$ in $L^{\infty}(0,T; H^1) \cap L^2(0,T; H^2)$. 
Now we prove the uniqueness of strong solutions. Suppose $(q_1, u_1)$ and $(q_2, u_2)$ are strong solutions of \eqref{intro1}--\eqref{qinit} with same initial data. Let $q = q_1 - q_2, u = u_1 - u_2$ and $p = p_1 - p_2$. Then $q$ satisfies 
\be \la{uniqq}
\pa_t q + \l q = - u_1 \cdot \na q - u \cdot \na q_2
\ee and $u$ satisfies
\be \la{uniqu}
\pa_t u - \Delta u + \na p = - qRq_1 - q_2Rq - u_1 \cdot \na u - u \cdot \na u_2.
\ee 
We take the $L^2$ inner product of \eqref{uniqq} with $\l^{-1}q$ and the $L^2$ inner product of \eqref{uniqu} with $u$. We add the resulting energy equalities. We have a cancellation 
\be \la{uinitial}
- \int_{\RR^2} (u \cdot \na q_2) \l^{-1} q dx - \int_{\RR^2} (q_2Rq) \cdot u dx = 0 
\ee obtained from integration by parts. In view of the Ladyzhenskaya's interpolation inequality, we estimate
\be 
\left|\int_{\RR^2} (qRq_1) \cdot u dx \right|
\le C\|q\|_{L^2}\|q_1\|_{L^4}\|u\|_{L^2}^{\fr{1}{2}} \|\na u\|_{L^2}^{\fr{1}{2}}
\le \fr{1}{4} \|\na u\|_{L^2}^2 + \fr{1}{4} \|q\|_{L^2}^2 + C\|q_1\|_{L^4}^4 \|u\|_{L^2}^2
\ee and 
\be 
\left|\int_{\RR^2} (u \cdot \na u_2) \cdot u  dx \right|
\le \|u\|_{L^4}^2 \|\na u_2\|_{L^2}
\le \fr{1}{4} \|\na u\|_{L^2}^2 + C\|\na u_2\|_{L^2}^2 \|u\|_{L^2}^2.
\ee
Now we write
\be 
\int_{\RR^2} (u_1 \cdot \na q)  \l^{-1}q dx   
= \int_{\RR^2} \left(\l^{-\fr{1}{2}}(u_1 \cdot \na q) - u_1 \cdot \na \l^{-\fr{1}{2}}q \right) \l^{-\fr{1}{2}}q dx 
\ee via integration by parts, and we show below that 
\be\la{ufinal}
\left|\int_{\RR^2} \left(\l^{-\fr{1}{2}}(u_1 \cdot \na q) - u_1 \cdot \na \l^{-\fr{1}{2}}q \right) \l^{-\fr{1}{2}}q dx  \right|
\le  C\|u_1\|_{H^2} \|q\|_{L^2}\|\l^{-\fr{1}{2}}q\|_{L^2}.
\ee
Putting \eqref{uinitial}-\eqref{ufinal} together, we obtain the energy inequality 
\be 
\fr{d}{dt} \left[\|\l^{-\fr{1}{2}}q\|_{L^2}^2 + \|u\|_{L^2}^2 \right]
\le C\left[\| u_1\|_{H^2}^2 + \|\na u_2\|_{L^2}^2 + \|q_1\|_{L^4}^4 \right] \left[\|\l^{-\fr{1}{2}}q\|_{L^2}^2 + \|u\|_{L^2}^2 \right]
\ee from which we obtain uniqueness. 
Finally, we show that the estimate \eqref{ufinal} holds by establishing the commutator estimate 
\be  \la{commutator}
\|\l^{-\fr{1}{2}}(u_1 \cdot \na q) - u_1 \cdot \na \l^{-\fr{1}{2}}q  \|_{L^2} 
\le C\| u_1\|_{H^2} \|q\|_{L^2}.
\ee 
Indeed, let $w\in L^2(\RR^2)$. By Parseval's identity, we have 
\be 
\int_{\RR^2} (\l^{-\fr{1}{2}} (u_1 \cdot \na q) - u_1 \cdot \na \l^{-\fr{1}{2}}q)(x) w (x) dx
=  \int_{\RR^2} \mathcal{F}(\l^{-\fr{1}{2}} (u_1 \cdot \na q) - u_1 \cdot \na \l^{-\fr{1}{2}}q)(\xi) \mathcal{F}w (\xi) d\xi.
\ee But 
\be 
\mathcal{F} (\l^{-\fr{1}{2}}(u_1 \cdot \na q))(\xi) 
= \int_{\RR^2} |\xi|^{-\fr{1}{2}} (\xi \cdot \mathcal{F}u_1 (\xi - y)) \mathcal{F}q (y) dy
\ee and 
\be 
\mathcal{F} (u_1 \cdot \na \l^{-\fr{1}{2}}q) (\xi)
= \int_{\RR^2} |y|^{-\fr{1}{2}} (\xi \cdot \mathcal{F}u_1 (\xi - y)) \mathcal{F}q (y) dy.
\ee 
Consequently,
\beg{align}
&\left|\int_{\RR^2} (\l^{-\fr{1}{2}} (u_1 \cdot \na q) - u_1 \cdot \na \l^{-\fr{1}{2}}q)(x) w (x) dx\right| \nonumber
\\&\le \int_{\RR^2} \int_{\RR^2} \min \left\{|\xi|, |y| \right\} \left||\xi|^{-\fr{1}{2}} - |y|^{-\fr{1}{2}}  \right| |\mathcal{F} u_1 (\xi - y) | |\mathcal{F} q(y)| |\mathcal{F} w(\xi) | dy d\xi
\end{align} where we used 
\be 
|\xi \cdot \mathcal{F}u_1 (\xi - y)| \le \min \left\{|\xi|, |y| \right\} |\mathcal{F} u_1 (\xi - y)|
\ee which holds due to the fact that the velocity is divergence-free. 
We note that 
\be 
\min \left\{|\xi|, |y| \right\} \left||\xi|^{-\fr{1}{2}} - |y|^{-\fr{1}{2}}  \right|
\le \frac{\min \left\{|\xi|, |y| \right\}}{|\xi|^{\fr{1}{2}}|y|^{\fr{1}{2}}} |\xi - y|^{\fr{1}{2}}
\le |\xi - y|^{\fr{1}{2}}
\ee for all $\xi, y \in \RR^2$. Therefore, 
\beg{align} 
\left|\int_{\RR^2} (\l^{-\fr{1}{2}} (u_1 \cdot \na q) - u_1 \cdot \na \l^{-\fr{1}{2}}q)(x) w (x) dx\right| 
&\le \||.|^{\fr{1}{2}} \mathcal{F} u_1(.) \|_{L^1} \|q\|_{L^2}\|w\|_{L^2} \nonumber 
\\&\le C\|u_1\|_{H^2} \|q\|_{L^2}\|w\|_{L^2}
\end{align} by H\"older's inequality and Young's convolution inequality. This gives \eqref{commutator} completing the proof of Theorem \ref{strong}.\\

{\bf{Acknowledgment.}} The research of M.I. was partially supported be NSF grant DMS 2204614.\\

{\bf{Data Availability Statement.}} The research does not have any associated data.\\

{\bf{Conflict of Interest.}} The authors declare that they have no conflict of interest.\\


\begin{thebibliography}{99}


\bibitem{AI} E.~Abdo, M.~Ignatova, \emph{Long time dynamics of a model of electroconvection}, Trans. Amer. Math. Soc.~{\bf{374}}, 5849--5875 (2021).


\bibitem{AGSS} C.~Amrouche, V.~Girault, M.E.~Schonbek, T.P.~Schonbek, \emph{Pointwise decay of solutions and of higher derivatives to Navier-Stokes equations}, SIAM J. Math. Anal.~{\bf{31}}, 740--753 (2000).

\bibitem{BN} C.~Bjorland, C.J.~Niche, \emph{On the decay of infinite energy solutions to the Navier-Stokes equations in the plane}, Physica D: Nonlinear Phenomena~{\bf{240}} (7), 670--674 (2011). 

\bibitem{BS} C.~Bjorland , M.E.~Schonbek, \emph{On questions of decay and existence for the viscous camassa-holm equations}, Ann. I. H. Poincar\'e-NA~{\bf{25}}, 907-–936 (2008).

\bibitem{C} A.P.~Calder\'on, \emph{Singular Integrals}, Bull. Amer. Math. Soc.~{\bf{72}}, 427--465 (1966).


\bibitem{CC} A.~C\'ordoba, D.~C\'ordoba, \emph{A maximum principle applied to quasi-geostrophic equations}, Comm. Math. Phys.~{\bf{249}}, 511--528 (2004).



\bibitem{ceiv}  P.~Constantin, T.~Elgindi, M.~Ignatova, V.~Vicol, \emph{On some electroconvection models}, Journal of Nonlinear Science~{\bf{27}}, 197--211 (2017).


 

\bibitem{CW} P.~Constantin, J.~Wu, \emph{Behavior of solutions of 2D quasi-geostrophic equations}, SIAM J. Math. Anal.~{\bf{30}}, 937--948 (1999).



\bibitem{DDMB} Z.A.~Daya, V.B.~Deyirmenjian, S.W.~Morris, J.R.~de Bruyn, \emph{Annular electroconvection with shear}, Phys. Rev. Lett.~{\bf{80}}, 964--967 (1998).


\bibitem{DL} B.~Dong, Y.~Li, \emph{Large time behavior to the system of incompressible non-newtonian fluds in $\R^2$}, J. Math. Anal. Appl.~{\bf{298}},667-–676 (2004).


\bibitem{KM} R.~Kajikiya, T.~Miyakawa, \emph{On $L^2$ decay of weak solutions of the Navier-Stokes equations in $\R^n$}, Math. Z.~{\bf{192}}, 135--148 (1986).

\bibitem{KU1} I.~Kukavica, \emph{On the weighted decay for solutions of the Navier-Stokes system}, Nonlinear Analysis: Theory, Methods and Applications~{\bf{70}} (6), 2466--2470 (2009). 


\bibitem{KU} I.~Kukavica, \emph{Space-time decay for Solutions of the Navier-Stokes equations}, Indiana Univ.
Math. J.~{\bf{50}}, 205--222(2001).


\bibitem{L} J.~Leray, \emph{Sur le mouvement d'un liquide visquex emplissant l'espace}, Acta Math.~{\bf{63}}, 193--248(1934).

\bibitem{NS} C.J.~Niche , M.E.~Schonbek , \emph{Decay of weak solutions to the 2d dissipative quasi-geostrophic equation}, Comm. Math. Phys.~{\bf{276}}, 93-–115 (2007).

\bibitem{OT} M.~Oliver, E.S.~Titi, \emph{Remark on the rate of decay of higher order derivatives for solutions to the Navier-Stokes equations in $\R^n$}, J. Funct. Anal.~{\bf{172}}, 1--18 (2000).

\bibitem{S} M.E.~Schonbek, \emph{$L^2$ decay for weak solutions of the Navier-Stokes equations}, Arch. Rational
Mech. Anal.~{\bf{88}}, 209--222 (1985). 



\bibitem{S1} M.E.~Schonbek, \emph{Uniform decay rated for parabolic conservations laws}, Nonlinear Analysis: Theory, Methods and Applications~{\bf{10}}, 943-–956 (1986).

\bibitem{SW} M.E.~Schonbek, M.~Wiegner, \emph{On the decay of higher-order norms of the solutions of Navier-Stokes equations}, Proc. Roy. Soc. Edinburgh Sect. A~{\bf{126}}, 677--685 (1996).

\bibitem{T} S.~Takahashi, \emph{A weigthed equation approach to decay rate estimates for the Navier-Stokes equations}, Nonlinear Anal.~{\bf{37}}, 751--789 (1999).

\bibitem{DMT} P.~Tsai, Z.~Daya, S.~Morris, \emph{Charge transport scaling in turbulent electroconvection}, Phys. Rev E~{\bf {72}}, 046311-1-12 (2005).

\bibitem{DDMT} P.~Tsai, Z.A.~Daya, V.B.~Deyirmenjian, S.W.~Morris, \emph{Direct numerical simulation of supercritical annular electroconvection}, Phys. Rev E~{\bf{76}}, 1--11 (2007). 
 
\bibitem{W1} M.~Wiegner, \emph{Decay results for weak solutions of the Navier-Stokes equations on $\R^n$}, J. London Math. Soc.~{\bf{35}}, 303--313 (1987).



\bibitem{Z} X.~Zhao, \emph{Asymptotic behavior of solutions to a new hall-MHD system}, Acta Applicandae Mathematicae~{\bf{157}}, 205--216 (2018). 
 
\bibitem{ZL} C.~Zhao, B.~Li, \emph{Time decay rate of weak solutions to the generalized MHD equations in $\R^2$}, Appl. Math. Comput.~{\bf{292}}, 1--8 (2017).   
 

\end{thebibliography}
\end{document}